\documentclass[a4paper,12pt]{article}
\usepackage{a4wide}
\usepackage{theorem}
\usepackage{amsmath}
\usepackage{array}
\usepackage{amssymb}
\usepackage{amsfonts}
\usepackage[english]{babel}
\usepackage{epsf}
\usepackage{epsfig}
\usepackage{graphicx}

\newtheorem{theo}{\indent Theorem \newline}[section]
\newtheorem{defi}[theo]{\indent Definition\newline}
\newtheorem{rem}[theo]{\noindent Remark}
\newtheorem{prop}[theo]{\indent Proposition\newline}
\newtheorem{lemma}[theo]{\indent Lemma\newline}

 \def\N{{\mathbb{N}}}
\def\Z{{\mathbb{Z}}}

\def\R{{\mathbb{R}}}

\usepackage{eufrak}

\newcommand{\ind}{\mathop{\rm ind}\nolimits}

\newcommand{\coker}{\mathop{\rm coker}\nolimits}

\newlength{\indentation}%
\makeatletter
\newcommand\@makefntextsans[1]{%
    \parindent 0em%
    \noindent%
    \hb@xt@0em{\hss}%
    #1}
\def\footnotetextsans{%
     \@ifnextchar [\@xfootnotenextsans%
       {\@footnotetextsans}}
\def\@xfootnotenextsans[#1]{%
  \begingroup%
     \csname c@\@mpfn\endcsname #1\relax%
  \endgroup%
  \@footnotetextsans}
\long\def\@footnotetextsans#1{\insert\footins{%
    \reset@font\footnotesize%
    \interlinepenalty\interfootnotelinepenalty%
    \splittopskip\footnotesep%
    \splitmaxdepth \dp\strutbox \floatingpenalty \@MM%
    \hsize\columnwidth \@parboxrestore%
    \color@begingroup%
      \@makefntextsans{%
        \rule\z@\footnotesep\ignorespaces#1\@finalstrut\strutbox}
    \color@endgroup}}
\makeatother

\begin{document}

\cleardoublepage
\title{Associahedral categories, particles and Morse functor}
\author{Jean-Yves Welschinger}
\maketitle

\makeatletter\renewcommand{\@makefnmark}{}\makeatother
\footnotetextsans{Keywords: Morse theory, categories.}
\footnotetextsans{AMS Classification : 53D40.
}

{\bf Abstract:}

Every smooth manifold contains particles which propagate. These form objects and morphisms of a category
equipped with a functor to the category of Abelian groups, turning this into a $0+1$ topological field theory.
We investigate the algebraic structure of this category, intimately related to the structure of Stasheff's polytops, 
introducing the notion of associahedral categories. An associahedral category is preadditive and close to being
strict monoidal. Finally, we interpret Morse-Witten theory as a contravariant functor, the Morse functor,  to the homotopy category
of bounded chain complexes of particles.

\section*{Introduction}

In a recent paper \cite{WelsFloer} dealing with Lagrangian Floer theory, a category of open strings was associated to
every symplectic manifold of dimension at least four, which comes equipped with a functor to the category of Abelian groups, 
turning this into a $1+1$ topological field theory. Lagrangian Floer theory was then interpreted as a contravariant functor, the
Floer functor,  to the homotopy category of bounded chain complexes of open strings.
The purpose of the present paper is to detail the Morse theoretic counterpart of this formalism. 
To every smooth manifold is associated a category of particles whose morphisms are trajectories of such particles. 
It comes with a functor, the functor Coefficients, whose target is the category of Abelian groups and which satisfies
the axioms of a $0+1$ topological field theory.
The category of particles has basically the same structure as the category of open strings, a structure which is
intimately related to Stasheff's polytopes, namely the associahedra and multiplihedra. We propose an algebraic
presentation of this structure and call such categories associahedral. These are small, preadditive and
almost strict monoidal categories. Actually, these categories come with morphisms of cardinality and indices to $(\Z , +)$ which 
we use to twist the functorial property of strict monoidal categories, see Definition \ref{defassociahedralcategory}.
Chain complexes in such categories satisfy some Leibnitz rule, chain maps are twisted morphisms with respect to the monoidal 
structure and homotopies satisfy some twisted Leibnitz rule. In the same way as in \cite{WelsFloer}, we interpret
the Morse data, given by a Morse function, or by a finite
collection of Morse functions in Fukaya's version, together with a generic Riemannian metric, as objects of a category of conductors. Then, we interpret Morse-Witten theory as a contravariant functor from the category of conductors to the homotopy category
of bounded chain complexes of particles. Composition with the functor Coefficients then provides the usual
chain complexes of free Abelian groups. As in the work \cite{WelsFloer}, this is possible thanks to the following observation:
a Morse function gives more than finitely many Morse trajectories between
critical points with index difference one. It gives finitely many such trajectories together with canonically oriented Fredholm
operators. These Fredholm operators are given by the first variation of Morse equation and are thus just connections along
the trajectories.

We first give an algebraic treatment of associahedral categories, then introduce the category of particles and
finally define the Morse functor. The latter is only defined for closed smooth manifolds.\\

{\bf Acknowledgements: }

I am grateful to the French Agence nationale de la recherche for its support. 

\tableofcontents

\section{Associahedral categories}

\subsection{Definition and the functor Coefficients}

Recall that a category ${\cal P} = \big( \text{Ob} ( {\cal P} ) , \text{Hom} ( {\cal P} ) \big)$ is said to be small
when $\text{Ob} ( {\cal P} ) $ and $ \text{Hom} ( {\cal P} ) $ are sets and preadditive whenever
the set $\text{Hom}_{P^+ , P^-} ( {\cal P} ) \subset \text{Hom} ( {\cal P} ) $ of morphisms $P^+ \to P^-$ has the structure of an
Abelian group for every objects $P^+ , P^- \in \text{Ob} ( {\cal P} ) $. We denote by
$s , t : \text{Hom} ( {\cal P} ) \to \text{Ob} ( {\cal P} )$ the source and target maps and by
$\emptyset_{P^+ , P^-} \in \text{Hom}_{P^+ , P^-} ( {\cal P} ) $ the unit element.
A category ${\cal P}$ is said to be strict monoidal iff it is equipped with a functor
$\otimes :  {\cal P} \times {\cal P} \to {\cal P} $ which is associative and has a unit element.
We denote by $(P , P') \in \text{Ob} ( {\cal P} ) \times \text{Ob} ( {\cal P} ) \mapsto
P \otimes P' \in \text{Ob} ( {\cal P} )$ and  $(\gamma , \gamma' ) \in \text{Hom}_{P^+ , P^-} ( {\cal P} ) 
\times \text{Hom}_{P'^+ , P'^-} ( {\cal P} ) \mapsto \gamma \otimes \gamma'  \in 
\text{Hom}_{P^+ \otimes P'^+ , P^- \otimes P'^-} ( {\cal P} ) $ these associative products and by 
$\emptyset \in  \text{Ob} ( {\cal P} ) $,
$\text{id}_\emptyset \in \text{Hom}_{\emptyset , \emptyset} ( {\cal P} ) $ their unit elements.
Functoriality means that $\text{id}_P \otimes \text{id}_{P'} = \text{id}_{P \otimes P'}$ for every 
$P , P' \in \text{Ob} ( {\cal P} )$ and $(\gamma_2 \otimes \gamma_2' ) \circ (\gamma_1 \otimes \gamma_1' )
= (\gamma_2 \circ \gamma_1 ) \otimes (\gamma_2'  \circ \gamma_1' ) $ whenever $s (\gamma_2) =
t(\gamma_1) $ and $s (\gamma_2') = t (\gamma_1') $. We are going to twist the latter property in Definition
\ref{defassociahedralcategory}.

\begin{defi}
\label{defassociahedralcategory}
An associahedral category ${\cal P}$ is a small preadditive category equipped with an associative 
product $\otimes :  {\cal P} \times {\cal P} \to {\cal P} $  having a unit element $\emptyset \in  \text{Ob} ( {\cal P} ) $
and distributive with respect to the preadditive structure such that:

1) The set $(\text{Ob} ( {\cal P} ) , \otimes)$ is a free monoid equipped with morphisms of
cardinality $q : (\text{Ob} ( {\cal P} ) , \otimes) \to (\N , +)$ and index $\mu_N : (\text{Ob} ( {\cal P} ) , \otimes) \to (\Z/N\Z , +)$,
where $N \in \N$ is given. The unit element $\emptyset$ is the only one whose cardinality vanishes whereas
elements of cardinality one, called elementary, generate $\text{Ob} ( {\cal P} )$.

2) The set $(\text{Hom} ( {\cal P} ) , \otimes)$ is a free monoid generated by morphisms with elementary targets, whose
unit element is $\text{id}_\emptyset \in \text{Hom}_{\emptyset , \emptyset} ( {\cal P} ) $. It is equipped with a morphism
of index $\mu : (\text{Hom} ( {\cal P} ) , \otimes) \to (\Z , +)$ additive with respect to composition and 
such that $\mu_N \circ t - \mu_N \circ s = \mu \mod (N)$,
where the source and target maps are morphisms $s,t : (\text{Hom} ( {\cal P} ) , \otimes) \to (\text{Ob} ( {\cal P} ) , \otimes)$.

3) For every $P , P' \in \text{Ob} ( {\cal P} )$, $\text{id}_P \otimes \text{id}_{P'} = \text{id}_{P \otimes P'}$. Moreover,
for every $\gamma_1 ,\gamma_2 , \gamma'_1 , \gamma'_2 \in \text{Hom} ( {\cal P} ) $ such that $s (\gamma_2) =
t(\gamma_1) $ and $s (\gamma_2') = t (\gamma_1') $, we have the relation
$(\gamma_2 \otimes \gamma_2' ) \circ (\gamma_1 \otimes \gamma_1' )
= (-1)^{\mu (\gamma_1) \mu (\gamma'_2)}  (\gamma_2 \circ \gamma_1 ) \otimes (\gamma_2'  \circ \gamma_1' ) $.
\end{defi}

\begin{rem}
\label{rem3'}
The third property of Definition
\ref{defassociahedralcategory} could have been replaced by the property

3') For every $P , P' \in \text{Ob} ( {\cal P} )$, $\text{id}_P \otimes \text{id}_{P'} = \text{id}_{P \otimes P'}$. Moreover,
for every $\gamma_1 ,\gamma_2 , \gamma'_1 , \gamma'_2 \in \text{Hom} ( {\cal P} ) $ such that $s (\gamma_2) =
t(\gamma_1) $ and $s (\gamma_2') = t (\gamma_1') $, we have the relation
$(-1)^{\mu (\gamma'_1) \mu (\gamma_2)}  (\gamma_2 \otimes \gamma_2' ) \circ (\gamma_1 \otimes \gamma_1' )
=  (\gamma_2 \circ \gamma_1 ) \otimes (\gamma_2'  \circ \gamma_1' ) $.
 
 If ${\cal P}$ is an associahedral category, its opposite category ${\cal P}^{op}$, that is the category having same
 objects and morphisms but with source and target maps exchanged, satisfies property 3'. Likewise, if ${\cal P}$ and ${\cal P}^*$
 are equipped with  a bijection $* : {\cal P} \to {\cal P}^*$ such that $(P_1 \otimes \dots \otimes P_q)^* = P_q^* \otimes \dots \otimes P_1^*$,
 where ${\cal P}$ is associahedral, $ {\cal P}^*$ satisfies properties $1$ and $2$ of Definition
\ref{defassociahedralcategory} and $*$ preserves $\mu$. Then, $ {\cal P}^*$ satisfies property $3'$.

\end{rem}

We define the cardinality of a morphism $\gamma \in \text{Hom}_{P^+ , P^-} ( {\cal P} ) $ of an associahedral category ${\cal P}$ 
to be the quantity $q (\gamma) =  q (P^+) - q(P^-)$. 
Morphisms with elementary target -which generate $(\text{Hom} ( {\cal P} ) , \otimes)$ - are called elementary. The category
of open strings of a symplectic manifold, introduced in \cite{WelsFloer}, is associahedral. The aim of this paper is to introduce
likewise an associahedral category of particles associated to any smooth manifold.

Denote by ${\cal A}b$ the Abelian category of Abelian groups.

\begin{defi}
\label{deffunctorcoefficients}

Let ${\cal P}$ be an associahedral category. The functor ${\cal C} : {\cal P} \to {\cal A}b$ defined by
$P \in \text{Ob} ( {\cal P} ) \mapsto \text{Hom}_{\emptyset , P} ( {\cal P} ) \in \text{Ob} ( {\cal A}b )$ and
$\gamma \in \text{Hom}_{P^+ , P^-} ( {\cal P} )  \mapsto \gamma \circ : \gamma'  \in  \text{Hom}_{\emptyset , P^+} ( {\cal P} )  
\mapsto \gamma \circ \gamma'  \in  \text{Hom}_{\emptyset , P^-} ( {\cal P} )  $ is called the functor Coefficients.

\end{defi}

In our examples, all Abelian groups are free.

\subsection{Chain complexes}

\begin{defi}
\label{defchaincomplex}
Let ${\cal P}$ be an associahedral category. A bounded chain complex $(\Lambda , \delta)$ of elements of ${\cal P}$
is a finite set $\Lambda$ of objects of ${\cal P}$ together with a matrix $\delta : \Lambda \to \Lambda$ of morphisms
of  ${\cal P}$ having the following three properties.

$A_1)$ For every $P , P' \in \text{Ob} ( {\cal P} )$, $P \otimes P' \in \Lambda \implies P , P' \in \Lambda.$

$A_2)$ For every member $\gamma$ of $\delta$, $\mu (\gamma) + q (\gamma) = 1$ and $\emptyset \notin \Lambda$.

$A_3)$ For every $P_1 \otimes P_2 \in  \Lambda$, the matrix of opposite morphisms $\delta^{op}$ satisfies the Leibniz rule
$$\delta^{op}_{P_1 \otimes P_2} =  (-1)^{q(P_2) \mu(\delta^{op}_{P_1} )} \delta^{op}_{P_1} \otimes id_{P_2} + 
(-1)^{q(P_1)} id_{P_1} \otimes \delta^{op}_{P_2} .$$

\end{defi}

Axiom $A_1$ implies that $\Lambda$ contains all the elementary components of its objects. Axioms $A_2$ implies
that the function $\mu_N - q$ defines a graduation modulo $N$ on $\Lambda$ for which $\delta$ is of degree one.
Axioms $A_1$ and $A_3$ imply that the differential $\delta$ is determined by its elementary components, so that for every
$P_1 \otimes \dots \otimes P_q \in \Lambda$, where $P_i$ is elementary for every $1 \leq i \leq q$, 
$$\delta^{op}_{P_1 \otimes \dots \otimes P_q } = \sum_{i=1}^q (-1)^{(q-i) \mu(\delta^{op}_{P_i} ) + i-1}
id_{P_1 \otimes \dots \otimes P_{i-1} } \otimes \delta^{op}_{P_i} \otimes id_{P_{i+1} \otimes \dots \otimes P_{q} }.$$
It is convenient to consider  the pair $(\{ \emptyset \} , \emptyset_{\emptyset , \emptyset})$ as a
chain complex as well, it will be called the empty chain complex.

\begin{lemma}
\label{lemmachain1}
Let $\gamma_1 : P_1^+ \to P_1^-$ and $\gamma_2 : P_2^+ \to P_2^-$  be two morphisms of the associahedral category
${\cal P}$ and $(\Lambda , \delta)$  be a bounded chain complex of elements of ${\cal P}$ such that
$P_1^+ \otimes P_2^+ , P_1^- \otimes P_2^+ , P_1^+ \otimes P_2^- $ and $P_1^- \otimes P_2^-$ are in $\Lambda$
and such that $\gamma_1$ , $\gamma_2$ are members of $\delta$. Then, the contribution of $\gamma_1$ , $\gamma_2$
to the restriction of $\delta \circ \delta : P_1^+ \otimes P_2^+ \to P_1^- \otimes P_2^-$ vanishes.
\end{lemma}

{\bf Proof:}

From Leibniz rule $A_3$, the contribution of $\gamma_1$ , $\gamma_2$ to $\delta_{P_1^+ \otimes P_2^+}$ equals
$(-1)^{q(P_2^+) \mu(\gamma_1 )} \gamma_1 \otimes id_{P_2^+} \oplus (-1)^{q(P_1^+)} id_{P_1^+}  \otimes
\gamma_2$ whereas their contribution to $\delta_{P_1^- \otimes P_2^+}$, $\delta_{P_1^+ \otimes P_2^-}$ equals
$(-1)^{q(P_1^-)}  id_{P_1^-}  \otimes \gamma_2$ and $(-1)^{q(P_2^-) \mu(\gamma_1 )} \gamma_1 \otimes id_{P_2^-}$
respectively. From the third property of Definition \ref{defassociahedralcategory}, their contribution to
$\delta \circ \delta_{P_1^+ \otimes P_2^+}$ thus equals
$(-1)^{q(P_2^-) \mu(\gamma_1 ) + q(P_1^+)}  \gamma_1\otimes \gamma_2 +
(-1)^{\mu(\gamma_1 ) \mu(\gamma_2 ) +  q(P_2^+) \mu(\gamma_1 ) + q(P_1^-)} \gamma_1\otimes \gamma_2$.
Axiom $A_2$ then provides the result. $\square$

\begin{lemma}
\label{lemmachain2}
Let $\gamma_1 : P_2^+ \to P_2^0$ and $\gamma_2 : P_2^0 \to P_2^-$  be two morphisms of the associahedral category
${\cal P}$ such that $\gamma_2 \circ \gamma_1 = \emptyset$. Let $(\Lambda , \delta)$  be a bounded chain complex of elements of ${\cal P}$ such that $P_1 \otimes P_2^+ \otimes P_3$, $P_1 \otimes P_2^0 \otimes P_3$ and
$P_1 \otimes P_2^- \otimes P_3$ are in $\Lambda$ and such that $\gamma_1$ , $\gamma_2$ are members of $\delta$,
where $P_1 , P_2$ are objects of ${\cal P}$. Then, the contribution of $\gamma_1$ , $\gamma_2$
to the restriction of $\delta \circ \delta :  P_1 \otimes P_2^+ \otimes P_3 \to P_1 \otimes P_2^- \otimes P_3$ vanishes.
\end{lemma}

{\bf Proof:}

From Leibniz rule $A_3$, the contribution of $\gamma_1$ to $\delta_{P_1 \otimes P_2^+ \otimes P_3}$ equals
 $(-1)^{q(P_1) + q(P_3) \mu(\gamma_1 ) } id_{P_1} \otimes \gamma_1 \otimes  id_{P_3}$ while the contribution of
  $\gamma_2$ to $\delta_{P_1 \otimes P_2^0 \otimes P_3}$ equals 
  $(-1)^{q(P_1) + q(P_3) \mu(\gamma_2 ) } id_{P_1} \otimes \gamma_2 \otimes  id_{P_3}$. Hence, their contribution to
  the composition $\delta \circ \delta_{P_1 \otimes P_2^+ \otimes P_3}$ equals
  $(-1)^{q(P_3) \mu(\gamma_2 \circ \gamma_1 ) } id_{P_1} \otimes (\gamma_2 \circ \gamma_1) \otimes  id_{P_3}
  = \emptyset$.   $\square$

\subsection{Chain maps}

\begin{defi}
\label{defchainmaps}
Let ${\cal P}$ be an associahedral category. A chain map (or morphism) $H : (\Lambda , \delta) \to (\Lambda' , \delta')$ 
between the bounded chain complexes $(\Lambda , \delta)$  and $(\Lambda' , \delta')$ of ${\cal P}$
is a map satisfying the following two axioms: 

$B_1)$ For every member $\gamma$ of $H$ having non empty source, $\mu (\gamma) + q (\gamma) = 0$.
Moreover, $\delta'  \circ H = H \circ \delta$.

$B_2)$ For every $P_1 \otimes P_2 \in \Lambda' $, the opposite morphism  $H^{op}$ satisfies the relation
$H^{op}_{P_1 \otimes P_2} = (-1)^{q(P_2) \mu(H^{op}_{P_1} ) } H^{op}_{P_1} \otimes H^{op}_{P_2}.$

\end{defi}

Axiom $B_1$ implies that $H$ is of degree $0$ for the modulo $N$ graduation defined by the function $\mu_N - q$.
Axioms $A_1$ and $B_2$ imply that the morphism $H$ is determined by its elementary components, so that for every
$P_1 \otimes \dots \otimes P_q \in \Lambda'$, where $P_i$ is elementary for every $1 \leq i \leq q$, 
$$H^{op}_{P_1 \otimes \dots \otimes P_q} = (-1)^{\sum_{i=1}^q (q-i) \mu(H^{op}_{P_i})} H^{op}_{P_1} \otimes 
\dots \otimes H^{op}_{P_q}.$$

\begin{lemma}
\label{lemmacomposition}
Let ${\cal P}$ be an associahedral category and $H : (\Lambda , \delta) \to (\Lambda' , \delta')$,
$H' : (\Lambda' , \delta') \to (\Lambda'' , \delta'')$ be two chain maps given by Definition \ref{defchainmaps}.
Then, the composition $H' \circ H : (\Lambda , \delta) \to (\Lambda'' , \delta'')$ is a chain map.
\end{lemma}

{\bf Proof:}

This composition satisfies Axiom $B_1$. Moreover,

\begin{eqnarray*}
(H' \circ H)^{op}_{P_1 \otimes P_2} & =&  H^{op} \circ H'^{op}_{P_1 \otimes P_2} \\
& =& (-1)^{q(P_2) \mu(H'^{op}_{P_1} ) } (-1)^{(q(P_2) + q(H'^{op}_{P_2} ) ) \mu (H^{op}_{t \circ H'^{op}_{P_1}})}
(H^{op}_{t \circ H'^{op}_{P_1}} \otimes H^{op}_{t \circ H'^{op}_{P_2}}) \circ (H'^{op}_{P_1 } \otimes H'^{op}_{P_2}) \\
&&\quad \text{ from } B_2\\
& =& (-1)^{q(P_2) \mu((H' \circ H)^{op}_{P_1} ) } (H' \circ H)^{op}_{P_1} \otimes (H' \circ H)^{op}_{P_2},
\end{eqnarray*}
from $B_1$ and  the third property of Definition \ref{defassociahedralcategory}, so that this 
composition satisfies Axiom $B_2$. $\square$

\begin{lemma}
\label{lemmatwistedLeibniz}
Let $H : (\Lambda , \delta) \to (\Lambda' , \delta')$ be a degree $0$ map satisfying Axiom $B_2$ between
the bounded chain complexes $(\Lambda , \delta)$  and $(\Lambda' , \delta')$ of the associahedral category ${\cal P}$.
Then, for every $P_1 \otimes P_2 \in \Lambda' $, the morphisms $K = \delta'  \circ H$ or $ H \circ \delta$ satisfy
the following twisted Leibniz rule:
$K^{op}_{P_1 \otimes P_2} = (-1)^{q(P_2) \mu(K^{op}_{P_1} ) } K^{op}_{P_1} \otimes H^{op}_{P_2}
+ (-1)^{q(P_1) + (q(P_2) +1) \mu(H^{op}_{P_1} ) } H^{op}_{P_1} \otimes K^{op}_{P_2}$.
\end{lemma}

{\bf Proof:}

Let $P_1 \otimes P_2 \in \Lambda' $, we have
\begin{eqnarray*}
(H \circ \delta)^{op}_{P_1 \otimes P_2} & =& \delta^{op} \circ H^{op}_{P_1 \otimes P_2} \\
& =& (-1)^{q(P_2) \mu(H^{op}_{P_1} ) } \delta^{op} \circ (H^{op}_{P_1} \otimes H^{op}_{P_2}) \quad \text{ from } B_2\\
& =& (-1)^{q(P_2) \mu(H^{op}_{P_1} ) } (-1)^{(q(P_2) + q(H^{op}_{P_2} ) ) \mu (\delta^{op}_{t \circ H^{op}_{P_1}})} 
(\delta^{op}_{t \circ H^{op}_{P_1}}  \otimes id_{t \circ H^{op}_{P_2}} ) \circ (H^{op}_{P_1} \otimes H^{op}_{P_2}) \\
&&+ (-1)^{q(P_2) \mu(H^{op}_{P_1} ) } (-1)^{q(P_1) + q(H^{op}_{P_1} ) }
( id_{t \circ H^{op}_{P_1}}  \otimes \delta^{op}_{t \circ H^{op}_{P_2}}) \circ (H^{op}_{P_1} \otimes H^{op}_{P_2})
\quad \text{ from } A_3\\
& =& (-1)^{q(P_2) \mu((H \circ \delta)^{op}_{P_1} ) }  (H \circ \delta)^{op}_{P_1} \otimes H^{op}_{P_2} +
(-1)^{q(P_1) + (q(P_2) +1) \mu(H^{op}_{P_1} ) } H^{op}_{P_1} \otimes (H \circ \delta)^{op}_{P_2}, 
\end{eqnarray*}
from $A_2$, $B_1$ and the commutation relation 
$(H_{Q_1} \otimes H_{Q_2}) \circ (\delta \otimes id_{Q_2}) = (-1)^{\mu (\delta) \mu (H_{Q_2})} (H \circ \delta)
\otimes H_{Q_2}$ deduced from the third property of Definition \ref{defassociahedralcategory}, see Remark \ref{rem3'}.
Likewise,

\begin{eqnarray*}
(\delta' \circ H)^{op}_{P_1 \otimes P_2} & =& H^{op} \circ \delta'^{op}_{P_1 \otimes P_2} \\
& =& (-1)^{q(P_2) \mu(\delta'^{op}_{P_1} ) } H^{op} \circ (\delta'^{op}_{P_1} \otimes id_{P_2} ) + 
(-1)^{q(P_1)} H^{op} \circ (id_{P_1} \otimes \delta'^{op}_{P_2}) \quad \text{ from } A_3\\
& =& (-1)^{q(P_2) \mu(\delta'^{op}_{P_1} ) } (-1)^{q(P_2) \mu(H^{op}_{t \circ \delta'^{op}_{P_1}} ) } 
(H^{op}_{t \circ \delta'^{op}_{P_1}} \otimes H^{op}_{P_2}) \circ (\delta'^{op}_{P_1} \otimes id_{P_2} ) + \\
&& (-1)^{q(P_1)} (-1)^{(q(P_2) + q(\delta'^{op}_{P_2})) \mu(H^{op}_{P_1})} (H^{op}_{P_1} \otimes 
H^{op}_{t \circ \delta'^{op}_{P_2}} ) \circ (id_{P_1} \otimes \delta'^{op}_{P_2} ) \quad \text{ from } B_2\\
& =& (-1)^{q(P_2)  \mu ((\delta' \circ H)^{op}_{P_1})} (\delta' \circ H)^{op}_{P_1} \otimes H^{op}_{P_2} +
(-1)^{q(P_1) + (q(P_2) +1) \mu(H^{op}_{P_1} ) } H^{op}_{P_1} \otimes (\delta' \circ H)^{op}_{P_2},
\end{eqnarray*}
from $A_2$, $B_1$ and the same commutation relation as before. $\square$\\

Note that the twisted Leibniz rule given by Lemma \ref{lemmatwistedLeibniz} coincides with $A_3$ when
$H = id$. From this rule follows that a map $H$ given by Lemma \ref{lemmatwistedLeibniz} is a chain map
provided the equality $\delta'  \circ H = H \circ \delta$ holds for its elementary members.

\subsection{Homotopies}

\begin{defi}
\label{defprimitivehomotopies}
Let ${\cal P}$ be an associahedral category. A primitive homotopy $K : (\Lambda , \delta) \to (\Lambda' , \delta')$ 
between the bounded chain complexes $(\Lambda , \delta)$  and $(\Lambda' , \delta')$ of ${\cal P}$
is a map satisfying the following three axioms: 

$C_1)$ For every $P_1 \otimes \dots \otimes P_q \in \Lambda' $, the elementary morphism
$(\delta'  \circ K + K \circ \delta)^{op}_{P_i}$ is non trivial for at most one $i \in \{ 1 , \dots , q \}$.

$C_2)$ For every member $\gamma$ of $K$ having non empty source, $\mu (\gamma) + q (\gamma) = -1$.

$C_3)$ For every $P_1 \otimes P_2 \in \Lambda' $, the opposite morphism  $K^{op}$ satisfies the 
following twisted Leibniz rule:
$$K^{op}_{P_1 \otimes P_2} = (-1)^{q(P_2) \mu(K^{op}_{P_1} ) } K^{op}_{P_1} \otimes H^{op}_{P_2}
+ (-1)^{q(P_1) + (q(P_2) +1) \mu(H^{op}_{P_1} ) } H^{op}_{P_1} \otimes K^{op}_{P_2},$$
where $H : (\Lambda , \delta) \to (\Lambda' , \delta')$ is a chain map given by Definition \ref{defchainmaps}
\end{defi}

A primitive homotopy is determined by its elementary components, so that for every
$P_1 \otimes \dots \otimes P_q \in \Lambda'$, where $P_i$ is elementary for every $1 \leq i \leq q$, 
$$K^{op}_{P_1 \otimes \dots \otimes P_q} = \sum_{i=1}^q (-1)^{i - 1 + \sum_{j=1}^q (q-j) \mu(L^{op}_{P_j})
+ \sum_{j=1}^{i-1} \mu(H^{op}_{P_j})} H^{op}_{P_1} \otimes  \dots \otimes H^{op}_{P_{i-1}} \otimes  K^{op}_{P_i} 
\otimes H^{op}_{P_{i+1}} \otimes \dots \otimes H^{op}_{P_q},$$
where $L^{op}_{P_j} = H^{op}_{P_j}$ if $j \neq i$ and $L^{op}_{P_i} = K^{op}_{P_i}$.

\begin{lemma}
\label{lemmahomotopies}
Let $K : (\Lambda , \delta) \to (\Lambda' , \delta')$ be a primitive homotopy between the bounded chain complexes 
$(\Lambda , \delta)$  and $(\Lambda' , \delta')$ of the associahedral category ${\cal P}$ and
$H : (\Lambda , \delta) \to (\Lambda' , \delta')$ be the associated chain map. Then, 
$H + \delta'  \circ K + K \circ \delta$ is a chain map.
\end{lemma}

{\bf Proof:}

Let $P_1 \otimes P_2 \in \Lambda' $, we have
\begin{eqnarray*}
(K \circ \delta)^{op}_{P_1 \otimes P_2} & =& \delta^{op} \circ K^{op}_{P_1 \otimes P_2} \\
& =& (-1)^{q(P_2) \mu(K^{op}_{P_1} ) } \delta^{op} \circ (K^{op}_{P_1} \otimes H^{op}_{P_2}) 
+ (-1)^{q(P_1) + (q(P_2) +1) \mu(H^{op}_{P_1} ) } \delta^{op} \circ (H^{op}_{P_1} \otimes K^{op}_{P_2})\\
& =& (-1)^{q(P_2) \mu(K^{op}_{P_1} ) } \big( (-1)^{(q(P_2) + q(H^{op}_{P_2} ) ) \mu (\delta^{op}_{t \circ K^{op}_{P_1}})} 
(\delta^{op}_{t \circ K^{op}_{P_1}}  \otimes id_{t \circ H^{op}_{P_2}} ) \circ (K^{op}_{P_1} \otimes H^{op}_{P_2}) \\
&&+  (-1)^{q(P_1) + q(K^{op}_{P_1} ) }
( id_{t \circ K^{op}_{P_1}}  \otimes \delta^{op}_{t \circ H^{op}_{P_2}}) \circ (K^{op}_{P_1} \otimes H^{op}_{P_2}) \big) + \\
&& (-1)^{q(P_1) + (q(P_2) +1) \mu(H^{op}_{P_1} ) } \big( (-1)^{(q(P_2) + q(K^{op}_{P_2} ) ) 
\mu (\delta^{op}_{t \circ H^{op}_{P_1}})} 
(\delta^{op}_{t \circ H^{op}_{P_1}}  \otimes id_{t \circ K^{op}_{P_2}} ) \circ (H^{op}_{P_1} \otimes K^{op}_{P_2}) \\
&&+  (-1)^{q(P_1) + q(H^{op}_{P_1} ) }
( id_{t \circ H^{op}_{P_1}}  \otimes \delta^{op}_{t \circ K^{op}_{P_2}}) \circ (H^{op}_{P_1} \otimes K^{op}_{P_2}) \big)
\quad \text{ from } A_3\\
& =& (-1)^{q(P_2) \mu((K \circ \delta)^{op}_{P_1} ) }  (K \circ \delta)^{op}_{P_1} \otimes H^{op}_{P_2} +
(-1)^{q(P_2) \mu(H^{op}_{P_1} ) } H^{op}_{P_1} \otimes (K \circ \delta)^{op}_{P_2} + \\
&& (-1)^{q(P_1) + (q(P_2) +1) \mu((H \circ \delta)^{op}_{P_1} ) } (H \circ \delta)^{op}_{P_1} \otimes K^{op}_{P_2} +\\
&&(-1)^{q(P_1) + 1 + (q(P_2) +1) \mu(K^{op}_{P_1} ) } K^{op}_{P_1} \otimes (H \circ \delta)^{op}_{P_2}, 
\end{eqnarray*}
from the third property of Definition \ref{defassociahedralcategory}, see Remark \ref{rem3'}.
Likewise,
\begin{eqnarray*}
(\delta' \circ K)^{op}_{P_1 \otimes P_2} & =& K^{op} \circ \delta'^{op}_{P_1 \otimes P_2} \\
& =& (-1)^{q(P_2) \mu(\delta'^{op}_{P_1} ) } K^{op} \circ (\delta'^{op}_{P_1} \otimes id_{P_2} ) + 
(-1)^{q(P_1)} K^{op} \circ (id_{P_1} \otimes \delta'^{op}_{P_2}) \quad \text{ from } A_3\\
& =& (-1)^{q(P_2) \mu(\delta'^{op}_{P_1} ) } \big( (-1)^{q(P_2) \mu(K^{op}_{t \circ \delta'^{op}_{P_1}} ) } 
(K^{op}_{t \circ \delta'^{op}_{P_1}} \otimes H^{op}_{P_2}) \circ (\delta'^{op}_{P_1} \otimes id_{P_2} ) + \\
&& (-1)^{q(P_1) + q(\delta'^{op}_{P_1}) + (q(P_2) + 1)  \mu(H^{op}_{t \circ \delta'^{op}_{P_1}} )}
(H^{op}_{t \circ \delta'^{op}_{P_1}} \otimes K^{op}_{P_2}) \circ (\delta'^{op}_{P_1} \otimes id_{P_2} ) \big) + \\
&& (-1)^{q(P_1)} \big( (-1)^{(q(P_2) + q(\delta'^{op}_{P_2})) \mu(K^{op}_{P_1})} (K^{op}_{P_1} \otimes 
H^{op}_{t \circ \delta'^{op}_{P_2}} ) \circ (id_{P_1} \otimes \delta'^{op}_{P_2} ) + \\
&& (-1)^{q(P_1) + (q(P_2) + q(\delta'^{op}_{P_2}) + 1) \mu(H^{op}_{P_1})} (H^{op}_{P_1} \otimes 
K^{op}_{t \circ \delta'^{op}_{P_2}} ) \circ (id_{P_1} \otimes \delta'^{op}_{P_2} ) \big)
\quad \text{ from } C_3\\
& =& (-1)^{q(P_2)  \mu ((\delta' \circ K)^{op}_{P_1})} (\delta' \circ K)^{op}_{P_1} \otimes H^{op}_{P_2} +
(-1)^{q(P_2) \mu(H^{op}_{P_1} ) } H^{op}_{P_1} \otimes (\delta' \circ K)^{op}_{P_2} + \\
&&(-1)^{q(P_1) + 1 + (q(P_2) +1) \mu((\delta' \circ H)^{op}_{P_1} ) } (\delta' \circ H)^{op}_{P_1} \otimes K^{op}_{P_2} +\\
&&(-1)^{q(P_1) + (q(P_2) +1) \mu(K^{op}_{P_1} ) } K^{op}_{P_1} \otimes (\delta' \circ H)^{op}_{P_2},
\end{eqnarray*}
from the third property of Definition \ref{defassociahedralcategory}. Summing up, we deduce from $B_1$ the
relation $( \delta'  \circ K + K \circ \delta)^{op}_{P_1 \otimes P_2} =
(-1)^{q(P_2)  \mu ((\delta' \circ K+ K \circ \delta)^{op}_{P_1})} (\delta' \circ K+ K \circ \delta)^{op}_{P_1} \otimes 
H^{op}_{P_2} + (-1)^{q(P_2) \mu(H^{op}_{P_1} ) } H^{op}_{P_1} \otimes (\delta' \circ K+ K \circ \delta)^{op}_{P_2}$.
From $C_1$, one of the two terms in the sum vanishes, so that $H + \delta'  \circ K + K \circ \delta$ satisfies $B_2$
since $H$ does. Likewise, $H + \delta'  \circ K + K \circ \delta$ satisfies $B_1$. $\square$

\begin{defi}
\label{defhomotopies}
Two chain maps $H, H' : (\Lambda , \delta) \to (\Lambda' , \delta')$ given by Definition \ref{defchainmaps} are said to be
primitively homotopic iff there exists a primitive homotopy $K : (\Lambda , \delta) \to (\Lambda' , \delta')$  associated
to $H$ such that $H' = H + \delta'  \circ K + K \circ \delta$. They are said to be homotopic iff there are
chain maps $H_0, \dots , H_k : (\Lambda , \delta) \to (\Lambda' , \delta')$ such that $H_0 = H$, $H_k = H'$
and for every $1 \leq i \leq k$, $H_{i-1}$ and $H_i$ are primitively homotopic.
\end{defi}

Let ${\cal P}$ be an associahedral category. We denote by $\text{Ob} ( K^b ({\cal P}) )$ the set of
chain complexes given by Definition \ref{defchaincomplex} and by $\text{Hom} ( K^b ({\cal P}) )$ the set
of chain maps given by Definition \ref{defchainmaps} modulo homotopies given by
Definition \ref{defhomotopies}. The category $K^b ({\cal P}) = \big( \text{Ob} ( K^b ({\cal P}) ) , \text{Hom} ( K^b ({\cal P}) )
\big)$ is the homotopic category of bounded chain complexes of ${\cal P}$. We denote with slight abuse by
$\emptyset \in \text{Ob} ( K^b ({\cal P}) )$ the empty chain complex.

\begin{defi}
Let $(\Lambda , \delta) \in \text{Ob} ( K^b ({\cal P}) )$. An augmentation of $(\Lambda , \delta) $ is a chain map
$\emptyset \to (\Lambda , \delta) $. A complex equipped with an augmentation is said to be an augmented complex.
\end{defi}

\section{Particles}

\subsection{Metric ribbon trees and Stasheff's associahedron}
\label{subsectassociahedron}

Following \cite{BoardV}, for every integer $l \geq 2$, we denote by $K_l$ the space of (connected) metric trees $T_l$ satisfying the
following three properties.

1) Each edge of $T_l$ has a lengh in $[0,2]$ and a tree with an edge of length $0$ is identified with the tree obtained by
contraction of this edge.

2) The tree $T_l$ has $l+1$ monovalent vertices, labeled $v_0 , \dots , v_l$ where $v_0$ is the root, and
has no bivalent vertex.
The edges adjacent to the $l+1$ monovalent vertices have length one. The tree is oriented from the root to the leaves.

3) For every edge $e$ of $T_l$, the partition of $\{ v_0 , \dots , v_l \} \cong \{ 0 , \dots , l \}$ induced 
by the connected components of $T_l \setminus \{ e \}$ only contains cyclically convex intervals of the
form $[i, j]$, $i \leq j$, or $[j, l] \cup [0,i]$, $i \leq j$.

This space $K_l$ has the structure of a $(l-2)$-dimensional convex polytope of the Euclidian space
isomorphic to the associahedron of Stasheff, see \S $1.4$ of \cite{BoardV}. 

\begin{defi}
\label{defstablemetrictrees}
The elements of $K_l$, $l \geq 2$, are called metric ribbon trees.
\end{defi}
We agree that $K_1$ is a singleton and that the corresponding tree $T_1$ is a compact interval of length one.
Likewise, we agree that $K_0$ is a singleton and that the corresponding tree $T_0$ is an open-closed interval of length one,
isometric to $[0,1)$, that is with only one vertex, the root, and only one edge. These trees $T_0$ and $T_1$
will also be called metric ribbon trees.

The dihedral group $D_{l+1}$ acts on $K_l$ by cyclic permutation of the labeling. We denote by 
$\rho_l : K_l \to K_l$ the ordered $l+1$ automorphism  induced by the cyclic permutation 
$(0 , \dots , l) \mapsto (1 , \dots , l , 0)$. We denote by $\sigma_l : K_l \to K_l$ the ordered two  automorphism  
induced by the reflection which fixes zero. It decomposes as the product of  transpositions
$(1, l) (2, l-1) \dots (\frac{1}{2} l , \frac{1}{2} l + 1)$ when $l$ is even and as
$(1, l) (2, l-1) \dots (\frac{1}{2} (l-1) , \frac{1}{2} (l + 3))$ when $l$ is odd, so that $\frac{1}{2} (l+1)$ is fixed as well
in this case.
Finally, when $l$ is odd, we denote by $\tau_l : K_l \to K_l$ the ordered two  automorphism  
induced by the reflection $(0, l) (1, l-1) \dots (\frac{1}{2} (l-1) , \frac{1}{2} (l + 1))$ which has no fixed point.

\begin{lemma}
\label{lemmaorientations}
For every $l \geq 2$, $\rho_l $ acts as $(-1)^l$ on the pair of orientations of $K_l$ whereas the involution  
$\sigma_l$ acts as
$(-1)^{\lfloor \frac{1}{2}(l-1) \rfloor}$ and, when $l$ is odd, $\tau_l$ acts as $(-1)^{\frac{1}{2}(l+1)}$.
\end{lemma}

{\bf Proof:}

\begin{figure}[ht]
\begin{center}
\includegraphics[scale=0.5]{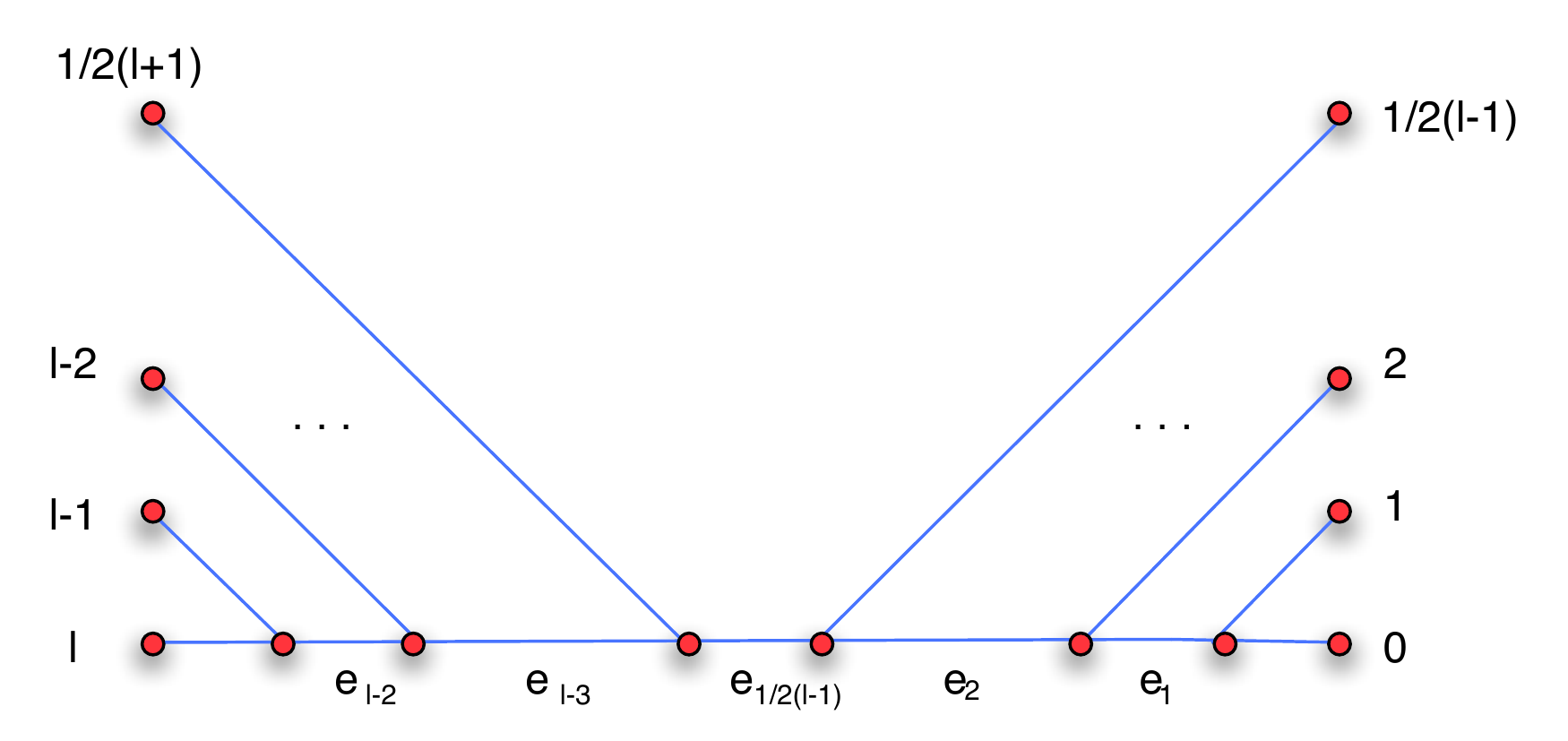}
\end{center}
\caption{}
\label{figurepeigne1}
\end{figure}

When  $l$ is odd, $\tau_l$ fixes the comb $P$ represented by Figure \ref{figurepeigne1} where every edge has
length one.

\begin{figure}[ht]
\begin{center}
\includegraphics[scale=0.5]{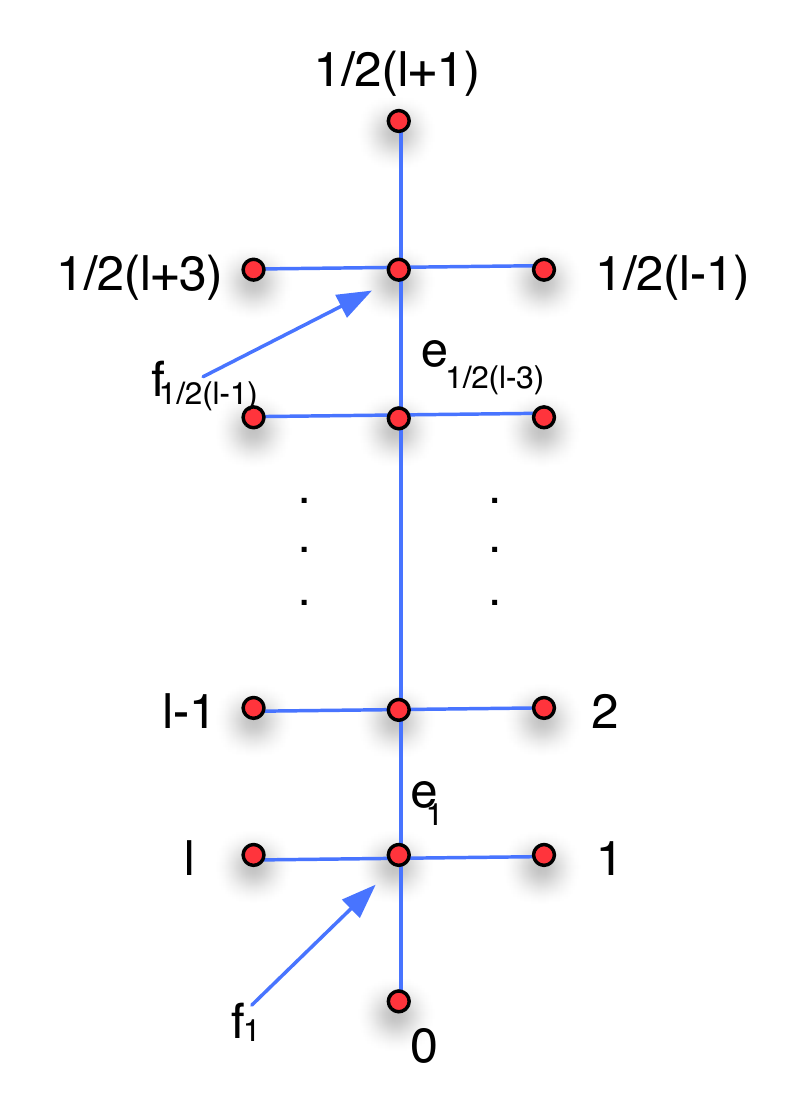}
\end{center}
\caption{}
\label{figurepeigne2}
\end{figure}

The tangent space $T_P K_l$ of $K_l$ at the comb $P$ is equipped with the basis 
$(e_1 , \dots , e_{l-2})$, where for $1 \leq i \leq l-2$, $e_i$ represents the stretch of the
edge represented by $e_i$ on Figure \ref{figurepeigne1}, fixing lengths of the other edges.
The differential map of $\tau_l$ at $P$ fixes the vector $e_{\frac{1}{2} (l-1)}$ and exchanges the 
pairs $(e_1 , e_{l-2}), \dots , (e_{\frac{1}{2} (l - 3)} , e_{\frac{1}{2} (l + 1)})$. Its action on
the pair of orientations of $T_P K_l$ is thus the same as the determinant of the antidiagonal matrix
of order $l-2$, that is $(-1)^{\frac{1}{2}(l+1)}$. Still when $l$ is odd, $\sigma_l$ fixes the tree $T$
represented by Figure \ref{figurepeigne2}.
The tangent space $T_T K_l$ of $K_l$ at this tree $T$ is equipped with the basis 
$(e_1 , \dots , e_{\frac{1}{2} (l - 3)} , f_1 , \dots , f_{\frac{1}{2} (l-1)})$, 
where for $1 \leq i \leq \frac{1}{2} (l - 3)$, $e_i$ represents the stretch of the
edge represented by $e_i$ on Figure \ref{figurepeigne2}, fixing lengths of the other edges,
 and for $1 \leq j \leq \frac{1}{2} (l - 1)$, $f_j$ represents the stretch of the length zero edge
 represented by the vertex of valence four denoted by $f_j$ on Figure \ref{figurepeigne2}.
The differential map of $\sigma_l$ at $T$ fixes all the vectors $e_i$ and reverses all the vectors
$f_j$, so that its action on the pair of orientations of $T_T K_l$ is $(-1)^{\frac{1}{2}(l-1)}$.
Finally, when $l$ is odd, $\rho_l = \sigma_l \circ \tau_l$, which finishes the proof of Lemma
\ref{lemmaorientations} in this case. Likewise, when $l$ is even, $\sigma_l$ acts on the 
the pair of orientations of $K_l$ as $(-1)^{\frac{1}{2}l-1}$ and $\rho_l $ decomposes as the
product of two elements conjugated to $\sigma_l$ so that it preserves orientations of $K_l$. $\square$\\

\begin{figure}[ht]
\begin{center}
\includegraphics[scale=0.5]{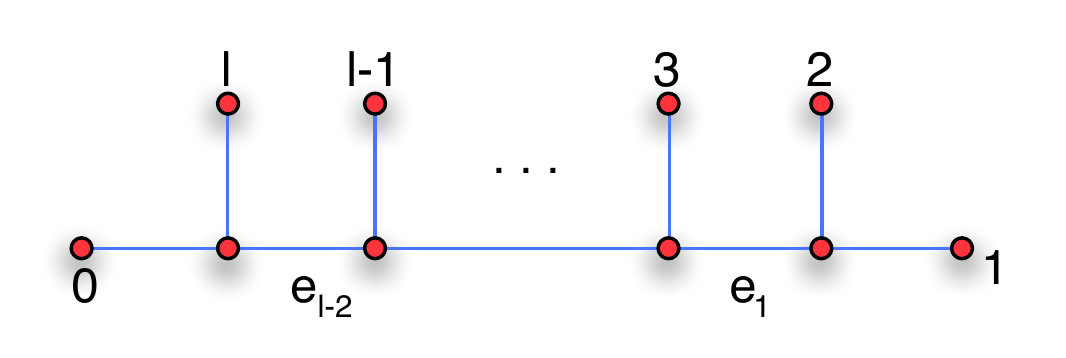}
\end{center}
\caption{}
\label{figurepeigne}
\end{figure}

From now on, for every $l \geq 2$, we denote by $P_l \in K_l$ the comb represented by Figure \ref{figurepeigne}.

The tangent space $T_{P_l} K_l$ of $K_l$ at the comb $P_l$ is equipped with the basis 
$(e_1 , \dots , e_{l-2})$, where for $1 \leq i \leq l-2$, $e_i$ represents the stretch of the
edge represented by $e_i$ on Figure \ref{figurepeigne}, fixing lengths of the other edges.
From now on, we equip $K_l$ with the orientation which turns this basis into a direct one.
Every face of $K_l$ canonically decomposes as a product of lower dimensional associahedra.
In particular, codimension one faces of $K_l$ decompose as products $K_{l_1} \times K_{l_2}$,
where $l_1 + l_2 = l+1$. More precisely, such a face writes
$F = \{ T = T' \cup T'' \in K_l \, \vert \, T' \in K_{l_1}  \, , \, T'' \in K_{l_2} \text{ and }  v'_i = v''_0 \}$,
where $i \in \{ 1 , \dots , l_1 \}$.
\begin{figure}[ht]
\begin{center}
\includegraphics[scale=0.5]{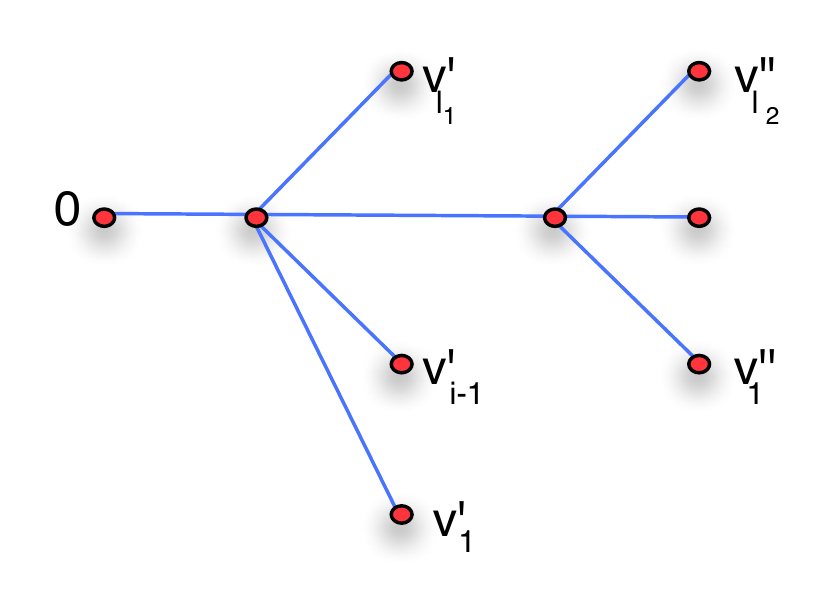}
\end{center}
\end{figure}
The only edge of length two  of the tree $T = T' \cup T'' \in \stackrel{\circ}{F}$ is the concatenation of the edge
adjacent to $v'_i$ in $T'$ with the one adjacent to the root $v''_0$ in $T''$. The labeled vertices $w_0 , \dots , w_l$
of $T$  satisfy $w_j = v'_j$ if $0 \leq j \leq i-1$, $w_j = v''_{j+1 - i}$ if $i \leq j \leq i + l_2 -1$ and
$w_j = v'_{j+1 - l_2}$ if $i + l_2 \leq j \leq l$. These facets inherit two orientations, one induced by $K_l$
and one induced by the product structure $K_{l_1} \times K_{l_2}$. Lemma \ref{lemmaassociahedron},
analogous to Lemma $2.9$ of \cite{WelsFloer}, compare these two orientations.

\begin{lemma}
 \label{lemmaassociahedron}
 Let $F = \{ T = T' \cup T'' \in K_l \, \vert \, T' \in K_{l_1}  \, , \, T'' \in K_{l_2} \text{ and }  v'_i = v''_0 \}$
 be a codimension one face of $K_l$, where $l_1, l_2 \geq 2$, $l_1 + l_2 = l+1$ and $1 \leq i \leq l_1$. 
 Then, the orientations of $F$ induced
 by $\partial K_l$ and by $K_{l_1} \times K_{l_2} $ coincide if and only if $l_1 l _2 + i(l_2 - 1)$ is odd.
 \end{lemma}

 {\bf Proof:}

Assume first that $i=1$. Then, the concatenation of the comb $P_{l_1} \in K_{l_1} $ with 
the comb $P_{l_2} \in K_{l_2} $ provides the comb $P_{l} \in K_{l} $. The infinitesimal
contraction of the length two edge $f$ of $P_l$ which joins $P_{l_1}$ and $P_{l_2}$
is an inward vector of $K_l$ at $P_l \in F$. The product orientation of $K_{l_1} \times K_{l_2} = F$
induces on $K_l$ the orientation which turns the basis
$(-f , e'_1 , \dots , e'_{l_1-2} , e''_1 , \dots , e''_{l_2-2})$ of $T_{P_l} K_l$ into a direct one.
The result follows in this case from the relation 
$f \wedge e'_1 \wedge \dots \wedge e'_{l_1-2} \wedge e''_1 \wedge \dots \wedge e''_{l_2-2} =
(-1)^{l_1 l_2 + l_2} e''_1 \wedge \dots \wedge e''_{l_2-2} \wedge f \wedge e'_1 \wedge \dots \wedge e'_{l_1-2}$.
Now the automorphism $\rho_l^{i-1}$ sends a facet for which $i > 1$ to a facet for which $i=1$.
This automorphism preserves the decomposition $K_{l_1} \times K_{l_2} $ of these facets, induces
the automorphism $\rho_{l_1}^{i-1}$ on the first factor $K_{l_1} $ and induces the identity on the second one.
The result thus follows from the first part and Lemma \ref{lemmaorientations}. $\square$

\begin{defi}
\label{defforest}
A forest $F_l$, $l \geq 0$, is an ordered finite union $T_{l_1} \cup \dots \cup T_{l_q} $
of metric ribbon trees $T_{l_i} \in K_{l_i}$, where $q \geq 1$, $l_1 + \dots + l_q = l$ and
$l_i \geq 0$ for every $1 \leq i \leq q$.
\end{defi}

For every forest $F_l = T_{l_1} \cup \dots \cup T_{l_q}$ given by Definition \ref{defforest}, we denote
by $v^-_1 , \dots , v^-_q$ the roots of $T_{l_1} , \dots , T_{l_q} $ respectively. Likewise, we denote
by $v^+_1 , \dots , v^+_l$ the leaves of $T_{l_1} \cup \dots \cup T_{l_q} $ in such a way that for
$1 \leq i \leq q$, $(v^-_i , v^+_{l_1 + \dots + l_{i-1} +1} , \dots , v^+_{l_1 + \dots + l_{i}})$ provides the
ordered set of labeled vertices of $T_{l_i}$.

Let $T_q \in K_q$, $q \geq 2$, be a metric ribbon tree and $F_l$ be a forest made of $q$ trees.
Let $T_q \cup F_l$ be the union of $T_q$ and $F_l$ where the leaves $v^+_1 , \dots , v^+_q$ of $T_q$
are identified with the roots $v^-_1 , \dots , v^-_q$ of $F_l$ respectively. If $F_l$ does not contain the trees
$T_0$ and $T_1$, this union $T_q \cup F_l$ is a metric ribbon tree denoted by $T_q \circ F_l$. For every copy of the tree
$T_1$ in $F_l$, the union $T_q \cup F_l$ contains an edge of length two adjacent to a leaf. We then divide
the metric of this edge by two to get the metric ribbon tree $T_q \circ F_l$. Likewise, for every copy of the tree
$T_0$ in $F_l$, the union $T_q \cup F_l$ contains a non compact edge of length two. We then remove 
(or contract) this free edge. In case $F_l$ only contains trees without leaf, we contract all such edges
except the first one for which we divide the metric by two. The tree $T_q \circ F_l$ coincides with $T_0$ in this case.
The composition $T_q \circ F_l$ is defined in the same way when $q=1$ or when $T_q$ is a forest instead of 
a tree.

\begin{defi}
\label{defcompositionforest}
Let $F_q$ and $F_l$ be forests such that $F_l$ is made of $q$ trees. The forest $F_q \circ F_l$ just defined above
is called the composition of $F_q$ and $F_l$.
\end{defi}

\subsection{The category of particles}

Let $M$ be a smooth manifold.

\begin{defi}
\label{defelementaryparticle}
An elementary particle $P$ of $M$ is a pair $(x , P_x)$, where $x$ is a point in $M$
and $P_x$ a linear subspace of $T_x M$. The point $x$ is called the based point of $P$
and denoted by $b(P)$.  The dimension $\mu (P)$ of $P_x$ is called the index of the elementary particle $P$.
\end{defi}

\begin{defi}
\label{defparticle}
A particle $P= P_1 \otimes \dots \otimes P_q$ of $M$ is an ordered finite collection of elementary
particles $P_i$ of $M$, $1 \leq i \leq q$. The integer $q$ is called the cardinality of the particle $P$
and denoted by $q (P)$. Its index $\mu (P)$ is the sum of the indices $\mu (P_i)$, $1 \leq i \leq q$.
The empty set $\emptyset$ is a particle whose cardinality and index vanish.
\end{defi}

We denote by $\text{Ob} ( {\cal P} (M) ) $ the set of particles of $M$. It is a free monoid for the product
$\otimes$, where we agree that $\emptyset$ is a unit element. The cardinality and index functions
define morphisms onto the monoid $(\N , +)$.

\begin{defi}
\label{defprimitivetrajectory}
A primitive trajectory from the particle $P^+= P^+_1 \otimes \dots \otimes P^+_l$ to the particle
$P^-= P^-_1 \otimes \dots \otimes P^-_q$ of $M$, $l \geq 0$,  $q \geq 1$, is a homotopy class with fixed ends of 
triples $\gamma = (u , F_l , \omega)$, such that:

1) $F_l =T_{l_1} \cup \dots \cup T_{l_q} $ is a forest with $q$ trees given by Definition \ref{defforest}.

2) $u : F_l \to M$ is a continuous map which sends the roots $v^-_1 , \dots , v^-_q$ of $T_{l_1} , \dots , T_{l_q} $  to the 
base points $b(P^-_1) , \dots , b(P^-_q) $ respectively
and their leaves $v^+_1 , \dots , v^+_l$ to the base points $b(P^+_1) , \dots , b(P^+_l) $ respectively.

3) $\omega$ is an orientation of the real line 
$ \big( \otimes_{i=1}^q \wedge^{\max} P^-_i \big) \otimes \big( \otimes_{j=1}^l \wedge^{\max} P^+_j \big)^*$.
\end{defi}

In particular, when $l=0$, such a primitive trajectory $\emptyset \to P^-$ reduces to an orientation of the linear space
$\oplus_{i=1}^q P^-_i $. If $\gamma = (u , F_l , \omega)$ is a primitive trajectory, we denote by $-\gamma$ the  primitive 
trajectory $(u , F_l , - \omega)$ obtained by reversing the orientation $\omega$. We agree that there are no primitive
trajectory with target the empty particle $\emptyset $ except the ones with empty source. There are two primitive trajectories
$\emptyset \to \emptyset $ denoted by $id_\emptyset $ and $-id_\emptyset $. For every particle $P$ of cardinality $l \geq 1$,
we denote by $id_P$ the primitive trajectory $(u , F_l , \omega) : P \to P$ for which $F_l$ is the union of $l$ copies of $T_1$,
$u$ is constant on every such tree and $\omega$ is the canonical orientation of the line $\wedge^{\max} P \otimes
(\wedge^{\max} P )^*$.

\begin{defi}
\label{deftrajectory}
A trajectory $\gamma = \oplus_{i \in I} \gamma_i : P^+ \to P^-$ from the particle 
$P^+= P^+_1 \otimes \dots \otimes P^+_l$ to the particle $P^-= P^-_1 \otimes \dots \otimes P^-_q$ of $M$ is a finite collection
of primitive trajectories $\gamma_i : P^+ \to P^-$, $i \in I$, modulo the relation 
$\gamma \cup (-\gamma) = \emptyset_{P^+ , P^-}$, where $\emptyset_{P^+ , P^-}$ is the empty trajectory.
The index $\mu(\gamma)$ of such a trajectory is the difference $\mu (P^-) - \mu (P^+)$
whereas its cardinality $q(\gamma)$ equals $q(P^+) - q(P^-)$. The trajectory is said to be elementary whenever
its target $P^-$ is an elementary particle.
\end{defi}

We denote by $\text{Hom} ( {\cal P} (M) ) $ the set of trajectories given by Definition \ref{deftrajectory}.
Trajectories between two particles come equipped with the operation $\oplus$ which turns this
set into a free Abelian group with unit the empty trajectory. If $\gamma' = (u' , F'_l , \omega') : P^+_1 \to P^-_1$
and $\gamma'' = (u'' , F''_l , \omega'') : P^+_2 \to P^-_2$ are two primitive trajectories, we denote by
$\gamma' \otimes \gamma'' : P^+_1 \otimes P^+_2  \to P^-_1\otimes P^-_2$ the primitive trajectory
$(u' \cup u'' , F'_l \cup F''_l , \omega' \otimes \omega'')$. Here, $\omega' \otimes \omega''$ is the orientation
induced on $ \big( \wedge^{\max} P^-_1 \otimes   \wedge^{\max} P^-_2 \big) \otimes \big(  \wedge^{\max} P^+_1 
 \otimes \wedge^{\max} P^+_2 \big)^*$ from $ \big( \wedge^{\max} P^-_1  \otimes (  \wedge^{\max} P^+_1)^* \big)
 \otimes \big( \wedge^{\max} P^-_2  \otimes (\wedge^{\max} P^+_2)^* \big)$ via the isomorphism
 $(  \wedge^{\max} P^+_1)^* \otimes   \wedge^{\max} P^-_2 \cong  \wedge^{\max} P^-_2 \otimes (  \wedge^{\max} P^+_1)^*$.
This product is extended to trajectories in such a way
that it is distributive with respect to the addition $\oplus$. Likewise, if $\gamma' = (u' , F'_l , \omega') : P^+ \to P^0$
and $\gamma'' = (u'' , F''_q , \omega'') : P^0 \to P^-$ are two primitive trajectories and $P^-$ is elementary,
we denote by $\gamma'' \circ \gamma' : P^+  \to P^-$ the primitive trajectory
$(u'' \cup u' , F''_q \circ F'_l , \omega'' \otimes \omega')$, where the composition of forests is the one given
by Definition \ref{defcompositionforest}. This operation $\circ$ is extended to primitive trajectories using the rule
$(\gamma''_1 \otimes \gamma''_2) \circ (\gamma'_1 \otimes \gamma'_2) = 
(-1)^{\mu (\gamma'_1) \mu (\gamma''_2)}  (\gamma''_1 \circ \gamma'_1) \otimes (\gamma''_2 \otimes  \gamma'_2)$.
It is then extended to trajectories in such a way
that it is distributive with respect to the addition $\oplus$.

\begin{prop}
\label{propparticles}
Let $M$ be a smooth manifold. The category of particles
${\cal P} (M) = \big( \text{Ob} ( {\cal P} (M) ) , \text{Hom} ( {\cal P} (M) ) \big)$ is associahedral.
\end{prop}

{\bf Proof:}

First of all ${\cal P} (M)$ is indeed a category. The composition of morphisms is given by the operation $\circ$, it
is indeed associative and has unit $id_P$ for every object $P$. This category is small, the preadditive structure
is given by the operation $\oplus$ and the product $\otimes : {\cal P} (M) \to {\cal P} (M)$ is associative and
distributive. Properties one and two of Definition \ref{defassociahedralcategory} are obviously satisfied with
$N = 0$. Property three also follows from our definitions. $\square$

\subsection{Coparticles and functor Coefficients}

The dual coparticle  of an elementary particle $P = (x , P_x)$ is the coparticle $P^* = (x , T_x M / P_x)$. Its index
satisfies the relation $\mu ( P) + \mu ( P^*) = \dim M$, so that $\mu ( \emptyset^*) = \dim M$. The dual
 coparticle of a particle $P = P_1 \otimes \dots \otimes P_l$ is $P^* = P^*_l \otimes \dots \otimes P^*_1$.
 We denote by $\text{Ob} ( {\cal P}^* (M) ) $ the set of coparticles of $M$.

\begin{defi}
\label{defprimitivecotrajectory}
A primitive cotrajectory from the coparticle $P^*= P^*_1 \otimes \dots \otimes P^*_q$ to the coparticle
$P'^*= P'^*_1 \otimes \dots \otimes P'^*_l$ of $M$, $l \geq 0$,  $q \geq 1$, is a homotopy class with fixed ends of 
triples $\gamma^* = (u^* , F_l , \omega^*)$, such that:

1) $F_l =T_{l_1} \cup \dots \cup T_{l_q} $ is a forest with $q$ trees given by Definition \ref{defforest}.

2) $u^* : F_l \to M$ is a continuous map which sends the roots $v^-_1 , \dots , v^-_q$ of $T_{l_1} , \dots , T_{l_q} $  to the 
base points $b(P^*_1) , \dots , b(P^*_q) $ respectively
and their leaves $v'^*_1 , \dots , v'^*_l$ to the base points $b(P'^*_1) , \dots , b(P'^*_l) $ respectively.

3) $\omega^*$ is an orientation of the real line 
$\big( \otimes_{j=1}^l \wedge^{\max} P'^*_j \big)  \otimes \big( \otimes_{i=1}^q \wedge^{\max} P^*_i \big)^*$.

A primitive cotrajectory from the coparticle $\emptyset^*$ to the coparticle
$P'^*= P'^*_1 \otimes \dots \otimes P'^*_l$ is an orientation of the real line $\big( \otimes_{j=1}^l \wedge^{\max} P'^*_j \big) $.
\end{defi}

Note that when $M$ is oriented, an orientation of the real line $\wedge^{\max} P^*$ is the same as an orientation
of $(\wedge^{\max} P)^*$, so that a  cotrajectory $\gamma^* : P^* \to P'^*$ is the same as a 
trajectory $\gamma : P' \to P$. A cotrajectory is a finite collection
of primitive cotrajectories $\gamma^*_i : P^* \to P'^*$, $i \in I$, modulo the relation 
$\gamma^* \cup (-\gamma^*) = \emptyset^*_{P^+ , P^-}$, compare Definition \ref{deftrajectory}.
We denote by $\text{Hom} ( {\cal P}^* (M) ) $ the set of cotrajectories of $M$. The pair
${\cal P}^* (M) = \big( \text{Ob} ( {\cal P}^* (M) ) ,\text{Hom} ( {\cal P}^* (M) ) \big)$ is the dual category of
${\cal P}  (M)$. 

Let ${\cal M}od_{\Z}$ be the category of free Abelian groups of finite type and ${\cal M}od^1_{\Z} \subset {\cal M}od_{\Z}$
be the subcategory of rank one such groups. The functor Coefficients provided by Definition \ref{deffunctorcoefficients}
give functors ${\cal C} : {\cal P} (M) \to {\cal M}od^1_{\Z}$ and ${\cal C}^* : {\cal P}^* (M) \to {\cal M}od^1_{\Z}$.
The two generators of the rank one free Abelian group $\Z_P$ (resp. $\Z_{P^*}$) associated to a particle $P$
(resp. coparticle $P^*$) are the two orientations of the real line $\wedge^{\max} P$ (resp $\wedge^{\max} P^*$).
When $M$ is oriented, the restrictions of ${\cal C}$ and ${\cal C}^*$ to elementary (co)particles give dual functors.

\section{Morse functor}

\subsection{Conductors}
\label{subsectionconductors}

\begin{defi}
\label{defconductors}
Let $M$ be a smooth manifold. A conductor of $M$ is a collection  $F = (f_0 , \dots , f_l ; g)$, where
$f_j - f_i : M \to \R$ are Morse functions with distinct critical points, $0 \leq i < j \leq l $, and
$g$ is a generic Riemannian metric  which is standard near every critical point of these functions.
\end{defi}

\begin{defi}
Let $F = (f_0 , \dots , f_l ; g)$ be a conductor of $M$. A subconductor of $F$ is a conductor of the form
$(f_{i_0} , \dots , f_{i_q} ; g)$, where $0 \leq i_0 < \dots < i_q \leq l $. We also say that $F$ is a refinement of
$(f_0 , \dots , f_q ; g)$ by $(f_{q+1} , \dots , f_l ; g)$, where $0 \leq q \leq l$.
\end{defi}

We denote by  $\text{Ob} ( {\cal C} (M) ) $ the set of conductors of $M$.

\begin{defi}
\label{defeffectivecontinuations}
An effective continuation $H$ from the conductor $F = (f_0 , \dots , f_q ; g)$ to the conductor
$\widetilde{F} = (\tilde{f}_0 , \dots , \tilde{f}_q ; \tilde{g})$ of $M$ is a homotopy $(F_t)_{t \in [0 , 1]}$
between the subconductor $F_0$ of $F$ associated to $I_H$ and the subconductor $F_1$ of $F$ 
associated to $\phi_H (I_H)$. Here, $I_H$ is a subset of $\{ 0 , \dots , q \}$ and $\phi_H :
I_H \to \{ 0 , \dots , l \}$ an increasing injective map.
\end{defi}

\begin{defi}
\label{defcontinuations}
A continuation $H : F \to \widetilde{F}$ from a conductor $F$ to a conductor
$\widetilde{F}$ of the smooth manifold $M$ is a homotopy class with fixed ends 
of effective continuations between these conductors.
\end{defi}

We denote by  $\text{Hom} ( {\cal C} (M) ) $ the set of continuations given by Definition \ref{defcontinuations}.
These notions of conductors and continuations are analogous to the one given in \cite{WelsFloer}.
Let $H : F \to \widetilde{F}$ be a continuation given by Definition \ref{defcontinuations}. The subconductor
of $\widetilde{F}$ associated to $\phi_H (I_H)$ is called the image of $H$ whereas the subconductor
of $F$ associated to $I_H$ is called the cokernel.
A sequence $F \stackrel{H}{\to} \widetilde{F} \stackrel{K}{\to} \hat{F}$ is called exact when the intersection
$\text{Im}{\phi_H} \cap I_K$ contains at most one element.

\begin{prop}
\label{propconductors}
Let $M$ be a smooth manifold. The pair ${\cal C} (M) = \big( \text{Ob} ( {\cal C} (M) ) , \text{Hom} ( {\cal C} (M) ) \big)$
has the structure of a small category equipped with the properties of sub-objects, refinements, exact sequences,
cokernel and image. 
\end{prop}

{\bf Proof:}

The composition of morphisms $\widetilde{H} \circ H$ satisfies, with the notations of Definition \ref{defeffectivecontinuations},
the relation $I_{\widetilde{H} \circ H} = I_H \cap \phi_H^{-1} (I_{\widetilde{H}} )$. It is obviously associative with
units.  $\square$

\subsection{Connections}

Let us denote by $I$ the closed interval $[a,b]$, where $-\infty \leq a < b \leq + \infty$. We denote by $H^{1,2} (I , \R^n)$
the Hilbert space of functions  of class $L^2$ with one derivative in $L^2$ in the sense of distributions.
Let $A : I \to M_n ( \R)$ be a bounded continuous map, we denote by $\partial_A = \frac{\partial}{\partial t} - A$ the
associated operator $H^{1,2} (I , \R^n) \to L^{2} (I , \R^n)$. We denote by $H^{1,2}_0 (I , \R^n) \subset H^{1,2} (I , \R^n)$
the subspace of functions  which vanish at the ends of $I$ (they are of class $C^{\frac{1}{2}}_{\text{loc}}$ from Sobolev's
Theorem). This is a closed subspace of vanishing index when $I = \R$, index $n$ when one of the extremities $a$ or $b$
is infinite and index $2n$ when $I$ is compact. For every invertible diagonalisable matrix $B \in GL_n ( \R)$, we denote
by $\mu^+ (B)$ (resp. $\mu^- (B)$) the number of its positive (resp. negative) eigenvalues, so that
$\mu^+ (B) + \mu^- (B) = n$.

\begin{lemma}
\label{lemmaR}
Let $A : \R \to M_n ( \R)$  be a constant path with value an invertible diagonalisable matrix. Then, the operator
$\partial_A : H^{1,2} (\R , \R^n) \to L^{2} (\R , \R^n)$ is an isomorphism.
\end{lemma}

{\bf Proof:}

Without loss of generality, we may assume that $n=1$. The operator then writes $\partial_A = \frac{\partial}{\partial t} - a$,
$a \in \R^*$, and is injective. The inverse operator writes $w \in  L^{2} (\R , \R^n) \mapsto 
\int_0^{+\infty} w (t-s) \exp (sa) ds \in H^{1,2} (\R , \R^n) $ when $a < 0$ and
$w \in  L^{2} (\R , \R^n) \mapsto - 
\int^0_{-\infty} w (t-s) \exp (sa) ds \in H^{1,2} (\R , \R^n) $ when $a > 0$. $\square$

\begin{lemma}
\label{lemmaI}
1) Let $A : \R^+ \to M_n ( \R)$  be a continuous map which converges to an invertible diagonalisable matrix $A_\infty$
at infinity. Then, the associated operator  $\partial_A : H^{1,2} (\R^+ , \R^n) \to L^{2} (\R^+ , \R^n)$ is Fredholm of index
$\mu^- (A_\infty)$.

2) Let $A : [a,b] \to M_n ( \R)$  be a continuous map, where $-\infty < a < b < + \infty$. Then, $\partial_A $ is Fredholm 
of index $n$.
\end{lemma}

{\bf Proof:}

In order to prove the second part, it suffices to prove that the restriction of $\partial_A$ to the subspace
$H^{1,2}_0 (I , \R^n)$ is Fredholm of index $-n$. This restriction is injective. Moreover, every
$v \in H^{1,2} (I , \R^n)$ satisfies the estimate $\Vert v \Vert_{H^{1,2}} = \Vert \frac{\partial}{\partial t} v \Vert_{L^{2}} 
+ \Vert v \Vert_{L^{2}} \leq \Vert \partial_A (v) \Vert_{L^{2}} + (\sup_{[a,b] } \Vert A \Vert + 1) \Vert v \Vert_{L^{2}} $.
From this elliptic estimate follows that the image of $\partial_A$ is closed. Indeed, let $w = \lim_{n \to +\infty} \partial_A (v_n)
\in L^{2} ([a,b] , \R^n)$ be a point in the closure of $\text{Im} (\partial_A)$. If $(v_n)_{n \in \N}$ is a bounded sequence
in $H^{1,2}_0 (I , \R^n)$, then it lies in a compact subset of $L^{2} ([a,b] , \R^n)$ from Rellich's theorem. From the above
elliptic estimate follows that a subsequence of $(v_n)_{n \in \N}$ is of Cauchy type and thus converges to
$v \in H^{1,2}_0 (I , \R^n)$. Hence, $w =\partial_A (v) $ lies in the image of $\partial_A$. If $(v_n)_{n \in \N}$ is not bounded,
we divide it by its norm and get in the same way a sequence converging to $v \in H^{1,2}_0 (I , \R^n)$ of norm one and
in the kernel of $\partial_A$, which is impossible. It remains to prove that the orthogonal complement of the image
of $\partial_A$ restricted to $H^{1,2}_0 (I , \R^n)$ is of dimension $n$. Let $w$ be a point in this complement. Then, for every
$v \in H^{1,2} (I , \R^n)$, $\int_a^b \left< \frac{\partial v}{\partial t} - A (v) , w \right> dt = 0$, so that
$\int_a^b \left< \frac{\partial v}{\partial t}, w \right> dt = \int_a^b \left< v , A^* (w) \right> dt $. As a consequence, the derivative
of $w$ in the sense of distributions is of class $L^2$, so that $w \in H^{1,2} ([a,b]  , \R^n)$. Moreover, an integration
by parts shows that $\frac{\partial w}{\partial t}  + A^* (w) = 0$, since $v$ vanishes at $a$ and $b$. The kernel of
this adjoint operator $\partial^*_A$ is $n$-dimensional made of the solutions of the linear system 
$\frac{\partial w}{\partial t}  = - A^* (w)$.

Let us now prove the first part of Lemma \ref{lemmaI}.
The restriction of $\partial_A$ to $H^{1,2}_0 (\R^+ , \R^n)$ is  injective. It suffices to prove that this restriction has closed image.
Indeed, the kernel of the adjoint operator is made of solutions of the linear system $\frac{\partial w}{\partial t}  = - A^* (w)$,
so that it is of finite dimension bounded by $n$. Hence, $\partial_A$ is Fredholm. There exists a continuous path 
$(A_s)_{s \in [0,1]}$ of bounded continuous maps $\R^+ \to M_n ( \R)$ converging to $A_\infty$
at infinity such that $A_0 = A$ and $A_1 \equiv A_\infty$ is constant. The kernel of $\partial^*_{A_1}$ is of dimension
$\mu^+ (A_\infty)$ so that $\ind (\partial_A) = \ind  (\partial_{A_1}) = n - \mu^+ (A_\infty) = \mu^- (A_\infty)$.
To prove that the restriction of $\partial_A$ to $H^{1,2}_0 (\R^+ , \R^n)$ has closed image we follow 
the scheme of Floer \cite{Floer}, compare \cite{schwMorse}. Let $M > 1$ and $\beta_M^+ : \R^+ \to \R$ be
a strictly increasing smooth function such that $\beta_M^+ (\tau) = 0 $ if $\tau \leq M$ and $\beta_M^+ (\tau) = 1 $ if $\tau \geq M+1$.
We set $\beta_M^0 = 1 - \beta_M^+ $. Let $B_M : \R \to M_n ( \R)$  be a continuous map such that $B_M (\tau) = A (\tau)$
if $\tau \geq M$ and $\sup_{\tau \in \R} \Vert B_M (\tau) - A_\infty \Vert  \leq \sup_{\tau \geq M} \Vert A (\tau) - A_\infty \Vert $.
From Lemma \ref{lemmaR}, there exists $M \gg 0$ such that $\partial_{B_M} : H^{1,2} (\R , \R^n) \to L^{2} (\R , \R^n)$ 
is an isomorphism and we choose such an $M$. There exist then constants $C_1^+ , C_2^+ > 0$ such that for every 
$v \in H^{1,2} (\R^+ , \R^n)$, 
$\Vert \beta_M^+ v \Vert_{H^{1,2}} \leq  C_1^+ \Vert \partial_{B_M} (\beta_M^+ v) \Vert_{L^{2}} = 
C_1^+ \Vert \partial_{A} (\beta_M^+ v) \Vert_{L^{2}} \leq C_1^+ \Vert \partial_{A} ( v) \Vert_{L^{2} (\R^+ , \R^n)} + 
C_2^+  \Vert v \Vert_{L^{2} ([0,M + 1] , \R^n)} $. Likewise, as before, there exist  constants $C_1^0 , C_2^0 > 0$ such that for every 
$v \in H^{1,2} (\R^+ , \R^n)$, 
$\Vert \beta_M^0 v \Vert_{H^{1,2}} \leq  \Vert \partial_{A} (\beta_M^0 v) \Vert_{L^{2} (\R^+ , \R^n)} + C_1^0
\Vert \beta_M^0 v \Vert_{L^{2} ([0,M + 1] , \R^n)} \leq \Vert \partial_{A} (v) \Vert_{L^{2} (\R^+ , \R^n)} + C_2^0
\Vert v \Vert_{L^{2} ([0,M+1] , \R^n)}$. Summing up, we deduce the elliptic estimate
$\Vert v \Vert_{H^{1,2}} \leq (C_1^+ +1) \Vert \partial_{A} ( v) \Vert_{L^{2} (\R^+ , \R^n)} + 
(C_2^+ + C_2^0)  \Vert v \Vert_{L^{2} ([0,M + 1] , \R^n)} $, which implies the result as before. $\square$\\

Now, let  $T_l \in K_l$ be a metric ribbon tree, where $l \geq 2$. We denote by ${\cal E}_{T_l}$ its set of edges and
by ${\cal V}_{T_l}$ its set of vertices. Every bounded edge of $T_l$ is by definition isometric to exactly one level of the function
$(x, y ) \in [0,1]^2 \mapsto xy \in [0, 1]$, the level zero if it is of length two and one if it is of length zero for example. 
We equip these edges with the measure $d\tau$ induced by the gradient of this function, so that the measure is infinite when they
are of length two. Likewise, the edge adjacent to the root $v_0$ (resp. leaves $v_1 , \dots , v_l$) is isometric to 
the $x$-axis (resp. $y$-axis) of $[0,1]^2$ and equipped with the induced infinite measure. This is called a gluing profile
in \cite{Hoferpolyfolds}. We denote by $H^{1,2} (T_l , \R^n)$
the Hilbert space of functions  of class $L^2$ with one derivative in $L^2$ in the sense of distributions for the measure $d\tau$ and 
by $H^{1,2}_0 (T_l , \R^n) \subset H^{1,2} (T_l , \R^n)$
the subspace of functions  which vanish at the vertices of $T_l$. The latter is canonically isomorphic to the product
$\Pi_{e \in {\cal E}_{T_l}} H^{1,2}_0 (e , \R^n)$. For every continuous map
$A : T_l \to M_n ( \R)$, we denote by $\partial_A = \frac{\partial}{\partial t} - A$ the
associated operator $H^{1,2} (T_l , \R^n) \to L^{2} (T_l , \R^n)$. 

\begin{lemma}
\label{lemmaT}
Let $T_l \in K_l$ be a metric ribbon tree and $A : T_l \to M_n ( \R)$ be a continuous map, where $l \geq 0$.
Assume that for $0 \leq i \leq l$, the value $A (v_i)$ of $A$ at the leaf $v_i$ is an invertible diagonalisable matrix $A_i$
and that the same holds at the middle of the length two edges of $T_l$. Then, the operator  $\partial_A : 
H^{1,2} (T_l , \R^n) \to L^{2} (T_l , \R^n)$ is Fredholm of index $\mu^+ (A_0) - \sum_{i=1}^l \mu^+ (A_i)$.
\end{lemma}

{\bf Proof:}

When $l=0$, the map $x \in \R^+ \mapsto -x \in \R^-$ conjugates the operator $\partial_A$ to the operator
$-\partial_{-A} : H^{1,2} (\R^+ , \R^n) \to L^{2} (\R^+ , \R^n)$ which is Fredholm of index $\mu^- (-A_0) =
\mu^+ (A_0)$ from Lemma \ref{lemmaI}. When $l \geq 2$, the restriction of $\partial_A$ to $H^{1,2}_0 (T_l , \R^n) $
is Fredholm of index $\mu^+ (A_0) + \sum_{i=1}^l \mu^- (A_i) - n  \# {\cal E}_{T_l}$ from Lemma \ref{lemmaI}.
Thus, $\partial_A$ is Fredholm of index $\mu^+ (A_0) + \sum_{i=1}^l \mu^- (A_i) - n (\# {\cal V}_{T_l} - \# {\cal E}_{T_l} - l - 1)
= \mu^+ (A_0) + \sum_{i=1}^l \mu^- (A_i) - n l $. The result follows along the same lines when $l=1$. $\square$\\

Let $\partial_A : 
H^{1,2} (T_l , \R^n) \to L^{2} (T_l , \R^n)$ be an operator given by Lemma \ref{lemmaT}.  Following 
the appendix of \cite{FloerHofer}, we denote by $\det (\partial_A)$ the real line $(\wedge^{\max} \ker \partial_A ) \otimes 
(\wedge^{\max} \coker \partial_A )^*$. For every $0 \leq i \leq l$, we denote by $P_i \subset \R^n$ the maximal linear subspace
invariant under $A_i$ on which $A_i$ has positive eigenvalues, so that $\dim (P_i) = \mu^+ (A_i)$.

\begin{lemma}
\label{lemmadet}
Let $\partial_A : H^{1,2} (T_l , \R^n) \to L^{2} (T_l , \R^n)$ be an operator given by Lemma \ref{lemmaT} and for every $0 \leq i \leq l$, 
$P_i \subset \R^n$ be the maximal linear subspace invariant under $A_i$ on which $A_i$ has positive eigenvalues.
Then, the real lines $\det (\partial_A)$ and $(\wedge^{\max} P_0) \otimes 
(\otimes_{i=1}^l \wedge^{\max} P_i )^*$ are canonically isomorphic.
\end{lemma}

{\bf Proof:}

Let us first assume that $l \geq 2$ and that $T_l$ does not contain any edge of length two.
For every $0 \leq i \leq l$, we denote by $e_i$ the length one edge adjacent to the leaf $v_i$. 
We deduce a decomposition $T_l = T'  \cup \cup_{i=0}^l e_i$, where $T'$ is a tree of finite measure.
We then deduce a short exact sequence of complexes from the sequence
$0 \to H^{1,2} (T_l , \R^n) \to H^{1,2} (T' , \R^n) \oplus \oplus_{i=0}^l  H^{1,2} (\overline{e}_i , \R^n)  \to \oplus_{i=0}^l \R^n_i \to 0$,
where the injective map is the restriction map to $T'$ and the closure $\overline{e}_i $ of $e_i$ and the surjective map
is the evaluation map at the intersections point $T' \cap \overline{e}_i $ of the difference between the two functions.
The isomorphism $\det (\partial_A \vert_{T'}) \otimes \otimes_{i=0}^l \det (\partial_A \vert_{\overline{e}_i}) \cong
\det (\partial_A) \otimes \det (\R^n_i)$ follows, see the appendix of \cite{FloerHofer}. Now the evaluation map at
$T' \cap \overline{e}_0$ provides an isomorphism between $\det (\partial_A \vert_{T'}) $  and $\det (\R^n_0)$ since from
Lemma \ref{lemmaI}, $\partial_A \vert_{T'}$ is of index $n$ and it has an $n$-dimensional kernel determined by the
initial condition at $T' \cap \overline{e}_0$. Likewise, evaluation at $T' \cap \overline{e}_0$ provides an isomorphism between 
$\det (\partial_A \vert_{\overline{e}_0}) $  and $\wedge^{\max} P_0$ whereas evaluation at $T' \cap \overline{e}_i$,
$1 \leq i \leq l$,  provides an isomorphism between 
$\det (\partial_A \vert_{\overline{e}_i}) \otimes \det (\R^n_i)^*$  and $(\wedge^{\max} P_i)^*$. Hence the result in this case. 
The result follows likewise from Lemma \ref{lemmaI} when $l=0$ or $1$ and by concatenation when $T_l$ contains edges of length two.
$\square$

\subsection{Gradient flow trajectories and Morse complexes}

Let $M$ be a closed smooth manifold and $F = (f_0 , \dots , f_l ; g) \in \text{Ob} ( {\cal C} (M) ) $ be a conductor of $M$.
For every $0 \leq i < j \leq l$ and every critical point $x$ of $f_j - f_i$, we denote by $P_x$ the elementary particle
based at the point $x$ whose linear space is the $+1$ eigenspace of the Hessian bilinear form 
$Hess_{f_{j} -  f_{i}} $ of $f_j - f_i$ for the metric $g$.
We set
$CM (f_i , f_j) = \oplus_{x \in \text{Crit} (f_i , f_j)} P_x$ and then for every $1 \leq q \leq l$,
$CM_q (F) = \oplus_{0 \leq i_0 < \dots < i_q \leq l} CM (f_{i_0} , f_{i_1}) \otimes \dots \otimes CM (f_{i_{q-1}} , f_{i_q})$.
Finally, we set $CM (F) = \oplus_{q=1}^l CM_q (F)$.

\begin{defi}
\label{defgradientflow}
Let $F = (f_0 , \dots , f_l ; g) \in \text{Ob} ( {\cal C} (M) ) $ be a conductor of the smooth manifold $M$. Let
$P_{x_1} \otimes \dots \otimes  P_{x_q}  \in CM (f_{i_0} , f_{i_1}) \otimes \dots \otimes CM (f_{i_{q-1}} , f_{i_q})$ and 
$P_{x_0}  \in CM (f_{i_0} , f_{i_q})$ be such that $\mu (P_{x_0}) - \mu(P_{x_1} \otimes \dots \otimes  P_{x_q}) = 2- q$.
A gradient flow trajectory $\gamma : P_{x_1} \otimes \dots \otimes  P_{x_q} \to P_{x_0} $ of $F$ is a 
primitive trajectory $(u, T , \omega)$ such that:

1) $T \in K_q$ is a metric ribbon tree with a root and $q$ leaves.

2) The restriction of $u : T \to M$ to every edge $e$ of $T$ satisfies for every $\tau \in e$ the equation
$\frac{\partial u}{\partial \tau} (\tau) = \nabla_g (f_{j_e} -  f_{i_e}) (u(\tau))$, where $d\tau$ is the measure element on $T$
and $[i_e+1, j_e] \cup ([j_e + 1, l] \cup [0,i_e ])$ is 
the partition of the leaves $\{ v_0 , \dots , v_q \} \cong \{ 0 , \dots , q \}$ into  cyclically convex intervals 
induced by $e$.

3) $ \omega$ is an orientation of the real line $\det (\nabla^\gamma)$, where $\nabla^\gamma : H^{1,2} (T , u^* TM) \to L^{2} (T , u^* TM)$
is the connection whose restriction to every edge $e$ of $T$ is the first variation
$v \mapsto \nabla_{\frac{\partial}{\partial \tau}} v  - Hess_{f_{j_e} -  f_{i_e}} (v)$. If $q > 1$, $\nabla^\gamma $ is injective
and has a $(q-2)$-dimensional cokernel canonically isomorphic to the tangent space $T_T K_q$. The orientation
 $ \omega$ is then the one induced from this isomorphism and the one chosen on $K_q$ in \S \ref{subsectassociahedron}.
 If $q=1$, $\nabla^\gamma $ is surjective and its one dimensional kernel is generated by $\frac{\partial u}{\partial \tau} $,
 $ \omega$ is then the orientation induced by  this vector.
\end{defi}
Recall that the tree $T \in K_q$ is oriented from the root to the leaves.
Note also that from Lemma \ref{lemmadet}, $\det (\nabla^\gamma)$ is canonically isomorphic to $(\wedge^{\max} P_{x_0}) \otimes 
(\otimes_{i=1}^l \wedge^{\max} P_{x_i} )^*$ so that a gradient flow trajectory given by Definition \ref{defgradientflow}
is indeed a primitive trajectory in the sense of Definition \ref{defprimitivetrajectory}.

From the genericness assumptions on $g$ and the compactness theorem in Morse theory, see \cite{schwMorse}, the conductor $F$ has only
finitely many gradient flow trajectories. We denote by $m_q : CM_q (F) \to CM_1 (F) $ the sum of all these trajectories.
Let then $\delta^{CM} : CM (F) \to CM (F)$ be the morphism induced from the opposite of these trajectories and axioms
$A_1 , A_2 , A_3$ of Definition \ref{defchaincomplex}. Hence, the restriction of 
$\delta^{CM}$ to $CM_q (F)$ equals
$\oplus_{l_2=1}^q (-1)^{ql_2} \oplus_{i=1}^{q_1} (-1)^{i(l_2 - 1)} id_{i-1} \otimes m_{l_2} \otimes id_{q_1 - i} \, ,$
where $q_1 + l_2 - 1 = q$, $1 \leq q \leq l$. The following theorem is due to Morse, Witten and Fukaya, see \cite{FukMorse},
\cite{FukOh}, \cite{schwMorse}.

\begin{theo}
\label{theomorsecomplex}
Let $M$ be a closed smooth manifold and $F \in \text{Ob} ( {\cal C} (M) ) $ be a conductor of $M$.
Then, $\delta^{CM} \circ \delta^{CM} = 0$, so that $(CM (F) , \delta^{CM}) \in \text{Ob} ( K^b ({\cal P} (M)))$.
\end{theo}

{\bf Proof:}

Let $P^+$ and $P^-$ in $CM (F)$ be particles of cardinalities $q^-$ and $q^+$ such that $\mu (P^-) - q^- = \mu(P^+) - q^+ + 2$, so that
the space ${\cal M}(P^+ ,P^- )$ of  pairs $(u,T)$ satisfying properties $1$ and $2$ of Definition \ref{defgradientflow}
is one-dimensional. From the compactness and glueing theorems in Morse theory, 
see \cite{schwMorse}, \cite{Hoferpolyfolds}, the union of primitive trajectories from $P^+$ to $P^-$ counted by 
$\delta^{CM} \circ \delta^{CM}$ 
 is in bijection with the boundary of this space. It suffices thus to prove that every such trajectory
 induces the outward normal orientation on ${\cal M}(P^+ ,P^- )$. Indeed, two boundary components of a same connected component
 of ${\cal M}(P^+ ,P^- )$ then provide the same trajectory but with opposite orientations, which implies the result.
 From Lemmas \ref{lemmachain1} and \ref{lemmachain2}, we can assume that $q^- = 1$, that is $P^-$ is elementary.
 Let $\gamma_1 \circ \gamma_2$ be such a trajectory counted by $\delta^{CM} \circ \delta^{CM}$, $P_0 = t(\gamma_2) = s(\gamma_1)$
 and $q^0$ be the cardinality of $P_0$. By definition, $\delta^{CM} \circ \delta^{CM}$ counts $\gamma_1 \circ \gamma_2$ with respect
 to the sign $1+ q^+ l_2 + i (l_2 - 1)$, where $ l_2 = q^+ + 1 - q^0 $, and with respect to the orientation induced by the product 
 $K_{q^0 } \times K_{l_2 }$ on the corresponding face of the associahedron $K_{q^+ }$.
 From Lemma \ref{lemmaassociahedron}, this orientation differs from the one induced by $K_{q^+ }$ exactly from the
 quantity $1+ q^0 l_2 + i (l_2 - 1) = 1+ q^+ l_2 + i (l_2 - 1) \mod(2)$. The result follows when $q_0 , l_2 \geq 2$.
 But the latter relation also remains valid when $q_0 $ or $ l_2 = 1$. Indeed, when $q_0 = 1$ (resp.  $ l_2 = 1$), 
 the associated connection $\nabla^{\gamma_1}$ (resp. $\nabla^{\gamma_2}$) is surjective and its one dimensional kernel
 is oriented such that it coincides with the outward (resp. inward) normal vector of ${\cal M}(P^+ ,P^- )$. This follows
 by derivation with respect to $R$ of the preglueing relation $u_1 \star_R u_2 (\tau) =  u_1 (\tau + R - R_0) \star_{R_0} u_2 
 (\tau - R + R_0)$ and from our
 convention to orient from target to source.  When $q_0 = 1$, the result is now obvious and when $ l_2 = 1$, it follows from the
 relation $k \wedge -\nu = (-1)^{q^0 - 1} \nu \wedge k$, where $k$ is an orientation of $K_{q^0 }$ and $\nu$ the outward normal
 vector of ${\cal M}(P^+ ,P^- )$. $\square$

\begin{defi}
\label{defmorsecomplex}
Let $M$ be a closed smooth manifold and $F \in \text{Ob} ( {\cal C} (M) ) $ be a conductor of $M$.
The  complex ${\cal M} (F) = (CM (F) , \delta^{CM}) \in \text{Ob} ( K^b ({\cal P} (M)))$ given by Theorem \ref{theomorsecomplex}
is called the Morse complex associated to $F$.
\end{defi}

\subsection{Morse continuations}

\subsubsection{Stasheff's multiplihedron}
\label{subsubsectmultiplihedron}

For every integer $l \geq 2$, denote by $J_l$ the space of 
metric ribbon trees with $l+1$ free edges given by Definition \ref{defstablemetrictrees} which are painted
in the sense of Definition $2.3.1$ of \cite{MauWood}, compare Definition $5.1$ of  \cite{Forc}. The trees are painted from the root to the leaves in such a way
that for every vertex $v$ of the tree, the amount of painting is the same on every subtree rooted at $v$
with respect to the measure $d\tau$. The space
 $J_l$ is a compactification of $ ] 0,1 [ \times \stackrel{\circ}{K}_l $ which has the structure of a $(l-1)$-dimensional convex 
 polytope of the Euclidian space
 isomorphic to Stasheff's multiplihedron, see \cite{Sta}, \cite{MauWood}, \cite{Forc} and references therin. 
 We agree that $J_0 = J_1 = [0,1] $. Forgetting the painting provides a map $for_l : J_l \to K_l$, we orient its fibers
 in the sense of propagation of the painting. This induces a product orientation on $\stackrel{\circ}{J}_l = ] 0,1 [ \times \stackrel{\circ}{K}_l $.
The codimension one faces
of $J_l$ different from $\{Ê0 \} \times K_l$ and $\{Ê1 \} \times K_l$ are of two different natures, see \cite{IwaMim}, \cite{Forc}.
The lower faces are canonically isomorphic to products $J_{l_1} \times K_{l_2}$, $l_1 + l_2 = l + 1$, they encode ribbon
trees having one edge of length two whose upper half is unpainted, see Definition $2.2$ of  \cite{Forc}. 
The upper faces are canonically isomorphic to products 
$K_q \times J_{l_1} \times  \dots \times  J_{l_q}$, $q \geq 1$, they encode ribbon
trees having $q$ edges of length two, whose lower half is painted, see Definition $2.3$ of  \cite{Forc}.
These faces inherit two orientations, one from $J_l$ and one from the product structure. The following Lemmas
\ref{lemmalowerface} and \ref{lemmaupperface}, analogous to Lemmas $3.4$ and $3.5$ of \cite{WelsFloer}, compare them.

\begin{lemma}
\label{lemmalowerface}
 Let $L = \{ T = T' \cup T'' \in J_l \, \vert \, T' \in J_{l_1}  \, , \, T'' \in K_{l_2} \text{ and }  v'_i = v''_0 \}$ be a lower facet
 of  the multiplihedron $J_l$, where $l_1 + l_2 = l+1$, $l_1, l_2 \geq 2$ and $1 \leq i \leq l_1$. Then, the orientations of $L$ induced
 by $J_l$ and by $J_{l_1} \times K_{l_2} $ coincide if and only if $l_1 l _2 + i(l_2 - 1)$ is even.
\end{lemma}

 {\bf Proof:}
 
 This result follows from Lemma \ref{lemmaassociahedron}. Indeed, at the infinitesimal level,
 $ ] 0,1 [ \times (K_{l_1} \times K_{l_2} ) = J_{l_1} \times K_{l_2} $ whereas $ ] 0,1 [ \times \partial K_l = - \partial J_l$. $\square$

\begin{lemma}
\label{lemmaupperface}
Let $U = \{   T  \cup  T^1 \cup \dots  \cup T^q \in J_l \, \vert \,  T    \in K_q \, , \, T^i \in J_{l_i} \text{ and } \, v_i = v_0^i
, 1 \leq i \leq q  \}$ be an upper facet
 of  the multiplihedron $J_l$, where $l_1 + \dots + l_q = l$, $q \geq 2$ and $l_i \geq 1$ for every $1 \leq i \leq q$. Then, 
 the orientations of $U$ induced
 by $J_l$ and $K_q \times J_{l_1} \times \dots \times J_{l_q}$ coincide if and only if $\sum_{i=1}^q (q-i)(l_i - 1)$ is even.
\end{lemma}

 {\bf Proof:}
 
 Let $T_q  \cup  P_{l_1} \cup \dots  \cup P_{l_q} \in U$, where $T_q \in K_q$ is a comb represented by Figure \ref{figurepeigne}
 which is completely painted and $P_{l_i} \in J_{l_i}$ is a comb represented by Figure \ref{figurepeigne} which is painted
 up to a level $s_i$. Let $T$ be an interior point of $J_l$ close to $T_q  \cup  P_{l_1} \cup \dots  \cup P_{l_q}$ so that
 its edges joining $T_q$ to $P_{l_1} , \dots  , P_{l_q}$  have length slightly less than two. The level of painting $s$ of $T$
 is an increasing function of $s_i$ and the length $f_i$ of the edge joining $T_q$ to $P_{l_i}$, $1 \leq i \leq q$.
 The interior of $U$ writes $\stackrel{\circ}{U} = \stackrel{\circ}{K}_q \times (]0,1[_1 \times  \stackrel{\circ}{K}_{l_1})
 \times \dots \times (]0,1[_q \times  \stackrel{\circ}{K}_{l_q})$ where each factor $]0,1[_i$ is equipped with the orientation
 $\frac{\partial}{\partial s_i}$, $1 \leq i \leq q$. Since for fixed $s$, increasing $f_i$ amounts to decrease $s_i$, this factor
 $]0,1[_i$ gets identified at $T$ to the length of the edge joining $T_q$ to $P_{l_i}$ but with the orientation 
 $- \frac{\partial}{\partial f_i}$ whereas $\frac{\partial}{\partial s}$ corresponds to the outward normal vector to $U$. We thus have 
 to compare the orientations of $\stackrel{\circ}{K}_l$ and $ \stackrel{\circ}{K}_q \times (]0,1[_1 \times  \stackrel{\circ}{K}_{l_1})
 \times \dots \times (]0,1[_q \times  \stackrel{\circ}{K}_{l_q})$ at $T$, where $]0,1[_i$ is equipped with the orientation 
 $- \frac{\partial}{\partial f_i}$. For this purpose, we proceed by finite induction. The glueing of $T_q $ with $P_{l_1}$
 gives the comb $P_{q + l_1 - 1}$. However, writing $(e_1 , \dots , e_{q + l_1 -3})$ the direct basis of
 $T_{P_{q + l_1 - 1}} K_{q + l_1 - 1}$ given by Figure \ref{figurepeigne}, we see that 
 $ \stackrel{\circ}{K}_q \times (]0,1[_1 \times  \stackrel{\circ}{K}_{l_1})$ comes equipped with the direct basis
 $( e_{ l_1} , \dots ,  e_{q + l_1 -3} ,   - e_{ l_1 -1}   ,   e_1 , \dots , e_{l_1 -2})$, so that the orientations of 
 $K_q  \times  J_{l_1}$ and $K_{q + l_1 - 1}$ differ from $(-1)^{(l_1 - 1)(q-1)}$. 
 
 Assume now that the orientations of
 $K_q \times J_{l_1} \times \dots \times J_{l_{i - 1}}$ and $K_{q + (l_1 - 1) + \dots + (l_{i - 1} - 1)}$ are compared and let
 us compare the orientations of $K_{q + (l_1 - 1) + \dots + (l_{i - 1} - 1)} \times J_{l_{i }}$ and $K_{q + (l_1 - 1) + \dots + (l_{i } - 1)}$.
 This general case differs from the first one because we now glue the comb $P_{q + (l_1 - 1) + \dots + (l_{i - 1} - 1)} $
 with $P_{l_{i }}$ at the vertex $v_I$ of $P_{q + (l_1 - 1) + \dots + (l_{i - 1} - 1)} $, where $I = 1 + l_1  + \dots + l_{i - 1} $, 
 so that the result is no more a comb. We can however reduce this case to the first one using the same trick as in Lemma
 \ref{lemmaassociahedron}, that is relabeling the vertices. Indeed, 
 the automorphism $\rho_{q + (l_1 - 1) + \dots + (l_{i } - 1)}^{I-1}$ of $K_{q + (l_1 - 1) + \dots + (l_{i } - 1)}$
 given by Lemma \ref{lemmaorientations} sends our facet onto a facet for which $I=1$.
This automorphism preserves the product structure of the facets and induces
the automorphism $\rho_{q + (l_1 - 1) + \dots + (l_{i - 1} - 1)}^{I-1}$ on $K_{q + (l_1 - 1) + \dots + (l_{i - 1} - 1)} $ 
and the identity on $J_{l_{i }}$. From Lemma \ref{lemmaorientations}, this automorphism
$\rho_{q + (l_1 - 1) + \dots + (l_{i } - 1)}^{I-1}$ contributes as $(-1)^{(l_{i } - 1)(I-1)}$ to the sign we are computing.
From the preceding case, we deduce that  the orientations of $K_{q + (l_1 - 1) + \dots + (l_{i - 1} - 1)} \times J_{l_{i }}$ and 
$K_{q + (l_1 - 1) + \dots + (l_{i } - 1)}$ differ from $(-1)^{(l_{i } - 1)(I-1 + q + (l_1 - 1) + \dots + (l_{i - 1} - 1) - 1)} = (-1)^{(l_{i } - 1)(q-i)}$. 
The result follows by summation. $\square$

\subsubsection{Morse continuations}

Let $H$ be an effective continuation from the conductor $F = (f_0 , \dots , f_l ; g) \in \text{Ob} ( {\cal C} (M) ) $ to the conductor
$\widetilde{F} = (\tilde{f}_0 , \dots , \tilde{f}_{\tilde{l}} ; g)$  of $M$, see \S \ref{subsectionconductors}. Restricting ourselves to 
the subconductor of $\widetilde{F}$ image of $H$, we may assume that $\tilde{l}=l$ and $I_H = \{ 0 , \dots , l \}$, so that $H$ is just
a homotopy $F^t = (f^t_0 , \dots , f^t_l ; g^t)$, $t \in [0 , 1]$, between $F^0 = F$ and $F^1 = \widetilde{F} $. Let  $1 \leq q \leq l$, 
$0 \leq i_0 < \dots < i_q \leq l$ and 
$P_{x_1} \otimes \dots \otimes  P_{x_q}  \in CM (\tilde{f}_{i_0} , \tilde{f}_{i_1}) \otimes \dots \otimes CM (\tilde{f}_{i_{q-1}} , \tilde{f}_{i_q})
\subset CM_q (\widetilde{F})$. Let 
$P_{x_0}  \in CM (f_{i_0} , f_{i_q})$ be such that $\mu (P_{x_0}) - \mu(P_{x_1} \otimes \dots \otimes  P_{x_q}) = 1- q$.
From the genericness assumptions on $g$ and the compactness theorem in Morse theory, see \cite{schwMorse}, 
\cite{Hoferpolyfolds}, there are only finitely many primitive trajectories 
$\gamma = (u, T , \omega) : P_{x_1} \otimes \dots \otimes  P_{x_q} \to P_{x_0} $ such that:

1) $T \in J_q$ is a metric ribbon tree with a root and $q$ leaves, painted up to a level $s$.

2) The restriction of $u : T \to M$ to every edge $e$ of $T$ satisfies for every $\tau \in e$ the equation

$$
\left\{
\begin{array}{rcl}
\frac{\partial u}{\partial \tau} (\tau) = \nabla_{\tilde{g}} (\tilde{f}_{j_e} -  \tilde{f}_{i_e}) (u(\tau)) &&\text{ if  $\tau \in e$ is unpainted},\\
\frac{\partial u}{\partial \tau} (\tau) = \nabla_g (f_{j_e} -  f_{i_e}) (u(\tau)) &&\text{ if  $\tau \in e$ is below the level $s-1$},\\
\frac{\partial u}{\partial \tau} (\tau) = \nabla_{g^t} (f^t_{j_e} -  f^t_{i_e}) (u(\tau)) &&\text{ if  $\tau \in e$ is at the level $s-t$},
\end{array}
\right.
$$
where $d\tau$ is the measure element on $T$, the levels $s-t$ are with respect to this measure
and $[i_e+1, j_e] \cup ([j_e + 1, l] \cup [0,i_e ])$ is 
the partition of the leaves $\{ v_0 , \dots , v_q \} \cong \{ 0 , \dots , q \}$ into  cyclically convex intervals 
induced by $e$.

3) $ \omega$ is an orientation of the real line $\det (\nabla^\gamma)$, where $\nabla^\gamma : H^{1,2} (T , u^* TM) \to L^{2} (T , u^* TM)$
is the connection whose restriction to every edge $e$ of $T$ is the first variation
$$
\left\{
\begin{array}{rcl}
v(\tau) \mapsto \nabla_{\frac{\partial}{\partial \tau}} v(\tau)  - Hess_{\tilde{f}_{j_e} -  \tilde{f}_{i_e}} (v(\tau)) &&
\text{ if  $\tau \in e$ is unpainted},\\
v(\tau) \mapsto \nabla_{\frac{\partial}{\partial \tau}} v(\tau)  - Hess_{f_{j_e} -  f_{i_e}} (v(\tau)) &&\text{ if  $\tau \in e$ is below the level $s-1$},\\
v(\tau) \mapsto \nabla_{\frac{\partial}{\partial \tau}} v(\tau)  - Hess_{f^t_{j_e} -  f^t_{i_e}} (v(\tau))  &&\text{ if  $\tau \in e$ is at the level $s-t$}.
\end{array}
\right.
$$
This connection $\nabla^\gamma $ is injective
and has a $(q-1)$-dimensional cokernel canonically isomorphic to the tangent space $T_T J_q$. The orientation
 $ \omega$ is then the one induced from this isomorphism and the one chosen on $J_q$ in \S \ref{subsubsectmultiplihedron}.
 From Lemma \ref{lemmadet} we know that $\det (\nabla^\gamma)$ is canonically isomorphic to $(\wedge^{\max} P_{x_0}) \otimes 
(\otimes_{i=1}^l \wedge^{\max} P_{x_i} )^*$ so that these triple are indeed primitive trajectories in the sense of Definition 
\ref{defprimitivetrajectory}.
  
 \begin{defi}
\label{defcontinuation}
These trajectories are the gradient trajectories of the continuation $H : F \to \widetilde{F} $.
\end{defi}

We denote by $H_q : CM_q (\widetilde{F}) \to CM_1 (F) $ the sum of all these trajectories of continuation.
Let then ${\cal M} (H) : CM (\widetilde{F}) \to CM (F)$ be the morphism induced by these elementary trajectories and axioms
$B_1 , B_2$ of Definition \ref{defchainmaps}. The restriction of 
${\cal M} (H)$ to $CM_{\tilde{q}} (\widetilde{F}) \to CM_{{q}} ({F})$ equals
$\oplus_{l_1 + \dots + l_q = \tilde{q}} (-1)^{\sum_{j=1}^q (q-j)(l_j - 1)} H_{l_1} \otimes 
\dots \otimes H_{l_q}$. In case $\tilde{l} > l$, we extend ${\cal M} (H)$ by zero outside of
the subcomplex associated to the subconductor of $\widetilde{F}$ image of $H$. The image of ${\cal M} (H)$ is the subconductor
of $CM (F)$ associated to the cokernel of $H$. 
The following theorem is due to Morse, Witten and Fukaya, see \cite{FukMorse}, \cite{schwMorse}.

\begin{theo}
\label{theocontinuation}
Let $M$ be a closed smooth manifold and $H$ be an effective continuation from the conductor $F$ to the conductor
$\widetilde{F}$ of $M$. Then, ${\cal M} (H) : CM (\widetilde{F}) \to CM (F)$ is a chain map.
\end{theo}

{\bf Proof:}

Let $P^+ \in CM (\widetilde{F})$ and $P^- \in CM (F)$ be particles of cardinalities $q^+$ and $q^-$ such that 
$\mu (P^-) - q^- = \mu(P^+) - q^+ + 1$, so that
the space ${\cal M}(P^+ ,P^- ; H )$ of  gradient trajectories of $H$ from $P^+$ to $P^-$ is one-dimensional. 
From the compactness and glueing theorems in Morse theory, 
see \cite{schwMorse}, \cite{Hoferpolyfolds}, the union of primitive trajectories from $P^+$ to $P^-$ counted by 
${\cal M} (H) \circ \delta^{CM} - \delta^{CM} \circ {\cal M} (H) $ 
 is in bijection with the boundary of this space. It suffices thus to prove that every such trajectory
 induces the outward normal orientation on ${\cal M}(P^+ ,P^- ; H)$. Indeed, two boundary components of a same connected component
 of ${\cal M}(P^+ ,P^- ; H )$ then provide the same trajectory but with opposite orientations, which implies the result.
 From Lemma \ref{lemmatwistedLeibniz}, we can assume that $q^- = 1$, that is $P^-$ is elementary.
 Let $\gamma_1 \circ \gamma_2$ be such a trajectory counted by $\delta^{CM} \circ {\cal M} (H)$, $P_0 = t(\gamma_2) = s(\gamma_1)$
 and $q^0$ be the cardinality of $P_0$. By definition, $\delta^{CM} \circ {\cal M} (H)$ counts $\gamma_1 \circ \gamma_2$ with respect
 to the sign $(-1)^{1+ \sum_{j=1}^{q^0} (q^0-j)(l_j - 1)}$, where $l_1 + \dots + l_{q^0}  = q^+  $, and with respect to the orientation 
 induced by the product $K_{q^0 } \times J_{l_1 } \times \dots \times J_{l_{q^0} }$
 on the corresponding upper face of the multiplihedron $J_{q^+ }$.
 From Lemma \ref{lemmaupperface}, this orientation differs from the one induced by $J_{q^+ }$ exactly from the
 quantity $(-1)^{\sum_{j=1}^{q^0} (q^0-j)(l_j - 1)}$, so that $\gamma_1 \circ \gamma_2$ is counted by $\delta^{CM} \circ {\cal M} (H)$
 with respect to the sign $-1$.
 This result remains valid when $q^0 = 1$, since in the neighborhood of the broken trajectory, 
 elements of ${\cal M}(P^+ ,P^- ; H )$ write $u_1 \star_R u_2 (\tau + R) =  u_1 (\tau + 2R - R_0) \star_{R_0} u_2 (\tau + R_0)$.
 Thus, $\lim_{R \to + \infty} \frac{\partial}{\partial \tau} u_1 \star_R u_2 (\tau + R) = 2 \frac{\partial u_1}{\partial \tau} $ and the vector
 $\frac{\partial u_1}{\partial \tau} $ which orients $m_1$ coincides with the outward normal orientation of ${\cal M}(P^+ ,P^- ; H )$.
  
 Likewise,  ${\cal M} (H) \circ \delta^{CM}$ counts  a trajectory $\gamma_1 \circ \gamma_2$ with respect
 to the sign $(-1)^{q^0 l_2 + i (l_2 - 1)}$,  where $q^0 + l_2  = q^+ + 1 $, and with respect to the orientation induced by the product 
 $J_{q^0 } \times K_{l_2 } $
 on the corresponding lower face of the multiplihedron $J_{q^+ }$.
 From Lemma \ref{lemmalowerface}, this orientation differs from the one induced by $J_{q^+ }$ exactly from the same
 quantity, so that $\gamma_1 \circ \gamma_2$ is counted by ${\cal M} (H) \circ \delta^{CM}$
 with respect to the sign $+1$.
 This result remains valid when $l_2 = 1$, since in the neighborhood of the broken trajectory, 
 elements of ${\cal M}(P^+ ,P^- ; H )$ write $u_1 \star_R u_2 (\tau - R) =  u_1 (\tau  - R_0) \star_{R_0} u_2 (\tau - 2R + R_0)$.
 Thus, $\lim_{R \to + \infty} \frac{\partial}{\partial \tau} u_1 \star_R u_2 (\tau - R) = -2 \frac{\partial u_2}{\partial \tau} $ and the vector
 $\frac{\partial u_2}{\partial \tau} $ which orients $m_1$ coincides with the inward normal orientation of ${\cal M}(P^+ ,P^- ; H )$.
 The result follows from the
 relation $k \wedge (-\nu) = (-1)^{q^0} \nu \wedge k$, where $k$ is an orientation of $J_{q^0 }$ and $\nu$ the inward normal
 vector of ${\cal M}(P^+ ,P^- ; H )$. $\square$

\subsection{Morse functor}

 \begin{theo}
\label{theofunctor}
Let $M$ be a closed smooth manifold. If $H^0, H^1$ are two effective continuations from the conductor $F$ to the conductor
$\widetilde{F}$ of $M$, then ${\cal M} (H^0)$ and $ {\cal M} (H^1) $ are chain homotopic maps
$CM (\widetilde{F}) \to CM (F)$. Likewise, if $H^0 : F^0 \to F^1$ and $H^1 : F^1 \to F^2$ are effective continuations,
then ${\cal M} (H^1 \circ H^0) = {\cal M} (H^0) \circ {\cal M} (H^1)$.
\end{theo}
Theorem \ref{theofunctor} means that ${\cal M} $ quotients out to a map 
$\text{Hom} (  {\cal C} (M) ) \to  \text{Hom} ( K^b ({\cal P} (M)))$ which together with the map
given by Definition \ref{defmorsecomplex} provide a contravariant functor
${\cal M} : {\cal C} (M)  \to K^b ({\cal P} (M))$.

\begin{defi}
\label{defmorsefunctor}
The contravariant functor ${\cal M} : {\cal C} (M)  \to K^b ({\cal P} (M))$ is called the Morse functor.
\end{defi}

{\bf Proof of Theorem \ref{theofunctor}:}

Let $(H^r)_{r \in [0,1]}$ be a generic homotopy between $H^0 $ and $H^1$. Equip $[0,1] \times J_l$ with the product orientation.
In the same way as we defined the morphisms ${\cal M} (H)$ in the previous paragraph, we define, 
for every $1 \leq q \leq l$, $0 \leq i_0 < \dots < i_q \leq l$ and $r \in [0,1]$, the morphisms
$K^r_q : CM_q (\widetilde{F}) \to CM_1 (F)$ as the sum of gradient trajectories of the continuation $H^r$ from a particle
$P^+ \in CM_q (\widetilde{F})$ to an elementary particle $P^- \in CM_1 (F)$ such that $\mu (P^-) - \mu(P^+)  = -q$.
These trajectories can only appear for a finite set of values of
$r \in ]0,1[$ since the homotopy $(H^r)_{r \in [0,1]}$ is generic. They
are oriented from the orientation just fixed on $[0,1] \times J_l$ since
the connection $\nabla^\gamma $ is injective
and has now a $q$-dimensional cokernel canonically isomorphic to the tangent space $T_T ([0,1] \times J_q)$.
Then, we denote by $K : CM (\widetilde{F}) \to CM (F)$ the morphism defined by
 $K = \oplus_{q=1}^l (-1)^q \oplus_{l_1 + \dots + l_q = \tilde{q}} (-1)^{\sum_{j=1}^{q} (q - j)(l_j - 1)}  
 \oplus_{i=1}^{q} (-1)^{\sum_{j=1}^{i-1} (l_j - 1)}
 \int_0^1  H^r_{l_1} \otimes \dots \otimes H^r_{l_{i-1}} \otimes K^r_{l_i} 
  \otimes H^r_{l_{i + 1}} \otimes \dots \otimes H^r_{l_{q}}$, where the integral is taken with respect to the counting measure.
 This morphism  $K$  is a homotopy in the sense of Definition \ref{defhomotopies} between ${\cal M} (H^0)$  and $ {\cal M} (H^1)$.
In other words, assuming that only one value $r_0$ of $r \in ]0,1[$ appear, $K$ is the morphism satisfying axioms 
$C_1 , C_2 , C_3$ of Definition \ref{defprimitivehomotopies}
whose opposite is given on elementary particles by $-K_q^{op}$, $1 \leq q \leq l$ and which satisfies the relation
${\cal M} (H^0) - {\cal M} (H^1) = \delta^{CM} \circ K + K \circ \delta^{CM} $.
Indeed, from the genericness of the homotopy $(H^r)_{r \in [0,1]}$, there is only one gradient trajectory 
$\gamma_0 : P_0^+ \to P_0^-$
of the continuation $H^{r_0}$. The elementary target particle $P_0^-$ belongs to a complex
$CM ({f}_{i_0 - 1} , {f}_{j_0})$. The morphism  $( K \circ \delta^{CM})^{op}$ is non trivial only on $P_0^-$
whereas $( \delta^{CM} \circ K)^{op}$ is non trivial only on complexes of the form $CM ({f}_{i - 1} , {f}_{j})$
with $i \leq i_0$ and $j \geq j_0$. For every particle $P_1 \otimes \dots \otimes P_q \in CM (F)$, at most one elementary
member belongs to $CM ({f}_{i_0 - 1} , {f}_{j_0})$ or such $CM ({f}_{i - 1} , {f}_{j})$, so that axiom
$C_1$ is satisfied. Axiom $C_2$ is by definition satisfied as well as axiom $C_3$ with $H = H^{r_0}$ or any $H^r$,
$r \in [0,1]$ actually. 

Let $P^+ \in CM (\widetilde{F})$ and $P^- \in CM (F)$ be particles of cardinalities $q^+$ and $q^-$ such that 
$\mu (P^-) - q^- = \mu(P^+) - q^+ $, so that
the space ${\cal M}(P^+ ,P^- ; K )$ of  gradient trajectories of $H^r$, $r \in [0,1]$, from $P^+$ to $P^-$ is one-dimensional. 
From the compactness and glueing theorems in Morse theory, 
see \cite{schwMorse}, \cite{Hoferpolyfolds}, the union of primitive trajectories from $P^+$ to $P^-$ counted by 
${\cal M} (H^1) - {\cal M} (H^0) + \delta^{CM} \circ K + K \circ \delta^{CM} $
 is in bijection with the boundary of this space. It suffices thus to prove that every such trajectory
 induces the outward normal orientation on ${\cal M}(P^+ ,P^- ; K)$. Indeed, two boundary components of a same connected component
 of ${\cal M}(P^+ ,P^- ; K )$ then provide the same trajectory but with opposite orientations, which implies the result.
 From Lemma \ref{lemmahomotopies}, we can assume that $q^- = 1$, that is $P^-$ is elementary.
 Let $\gamma_1 \circ \gamma_2$ be such a trajectory counted by $\delta^{CM} \circ K$, $P_0 = t(\gamma_2) = s(\gamma_1)$
 and $q^0$ be the cardinality of $P_0$. By definition, $\delta^{CM} \circ K$ counts $\gamma_1 \circ \gamma_2$ with respect
 to the sign $(-1)^{1+ q^0 + \sum_{j=1}^{q^0} (q^0-j)(l_j - 1) + \sum_{j=1}^{i-1} (l_j - 1)}$, where 
 $l_1 + \dots + l_{q^0}  = q^+  $, and with respect to the orientation induced by the product 
 $K_{q^0 } \times J_{l_1 } \times \dots \times J_{l_{i-1} } \times( [0,1] \times J_{l_i } )  \times J_{l_{i+1} } \times \dots \times J_{l_{q^0} }$
 on the corresponding face of  $[0,1] \times J_{q^+ }$.
 From Lemma \ref{lemmaupperface}, this orientation differs from the one induced by $[0,1] \times J_{q^+ }$ exactly from the
 quantity $(-1)^{\sum_{j=1}^{q^0} (q^0-j)(l_j - 1) + \sum_{j=1}^{i-1} (l_j - 1) + q^0 - 1}$, so that $\gamma_1 \circ \gamma_2$ is 
 counted by $\delta^{CM} \circ K$ with respect to the sign $+1$. The last $- 1$ comes from the commutation between
 the $[0,1] $-factor and the outward normal vector.
 This remains valid when $q^0 = 1$ for the same reason as in the proof of Theorem \ref{theocontinuation}.
  
 Likewise,  $K \circ \delta^{CM}$ counts  a trajectory $\gamma_1 \circ \gamma_2$ with respect
 to the sign $(-1)^{q^0 l_2 + i (l_2 - 1) + 1}$,  where $q^0 + l_2  = q^+ + 1 $, and with respect to the orientation induced by the product 
 $[0,1] \times J_{q^0 } \times K_{l_2 } $ on the corresponding face of $[0,1] \times J_{q^+ }$.
 From Lemma \ref{lemmalowerface}, this orientation differs from the one induced by $[0,1] \times J_{q^+ }$ exactly from the same
 quantity, so that $\gamma_1 \circ \gamma_2$ is counted by ${\cal M} (H) \circ \delta^{CM}$
 with respect to the sign $+1$. The last $+ 1$ comes once more from the commutation between
 the $[0,1] $-factor and the outward normal vector and 
 the result remains valid when $q^0 = 1$ for the same reason as in the proof of Theorem \ref{theocontinuation}.
 Hence the first part of Theorem \ref{theofunctor}.

Now, from Lemma \ref{lemmacomposition}, the composition ${\cal M} (H^0) \circ {\cal M} (H^1)$ is a chain map in the sense
of Defintion \ref{defchainmaps} and  from Lemma \ref{lemmatwistedLeibniz}, in order to prove the second part of 
Theorem \ref{theofunctor}, we just have to check that ${\cal M} (H^1 \circ H^0)$ and ${\cal M} (H^0) \circ {\cal M} (H^1)$
are homotopic on elementary members.
Let $\gamma_1 \circ \gamma_2$ be an elementary trajectory counted by ${\cal M} (H^0) \circ {\cal M} (H^1)$, 
$P_0 = t(\gamma_2) = s(\gamma_1)$ and $q^0$ be the cardinality of $P_0$. We write
$\gamma_2 = H^1_{l_1} \otimes \dots \otimes H^1_{l_{q_0}} $, where 
 $l_1 + \dots + l_{q^0} $ is the cardinality of the source of $\gamma_2$.
 By definition, ${\cal M} (H^0) \circ {\cal M} (H^1)$ counts $\gamma_1 \circ \gamma_2$ with respect
 to the sign $(-1)^{\sum_{j=1}^{q^0} (q^0-j)(l_j - 1) }$ and with respect to the orientation induced by the product 
 $J_{q^0 } \times J_{l_1 } \times \dots \times J_{l_{q_0}}$.
 This product is a face of the space $J^2_l$ which encodes painted metric ribbon trees with $l+1$  free edges among which 
 one root and with interior edges having lengths in $[0,2]$. But this time, these trees are painted with two different colors,
one which encodes the propagation of the homotopy $H^0$ and the other one which encodes the propagation of the homotopy $H^1$.
Thus, $J^2_l$ is the subset of the fiber product $J_l \times_{K_l} J_l$ made of pairs of painted trees with same underlying
tree and such that the first one is painted up to a lower level than the second one.
From the glueing theorem in Morse theory, see \cite{schwMorse}, \cite{Hoferpolyfolds},
this trajectory $\gamma_1 \circ \gamma_2$ is a boundary point of a one dimensional manifold whose underlying painted
tree belongs to the space $J^2_l$. The outward normal direction in this manifold is given by increasing the difference
between the two colors $1$ and $0$. Moreover, from the compactness theorem in Morse theory, see \cite{schwMorse}, 
\cite{Hoferpolyfolds}, for a value of this difference close to infinity, any such trajectory belongs to one such manifold.
Let us choose a generic value $s_1 - s_0$ close to infinity. We can define a morphism of continuation
$\widetilde{\cal M} (H^1 \circ H^0)$ in the same way as in Definition \ref{defcontinuation} by counting gradient trajectories
of $H^1 \circ H^0$ with difference of painting $s_1 - s_0$. Such elementary trajectories are counted with respect to the
sign $+1$ and the orientation of $J_{l_1 + \dots + l_{q^0} }$.  From Lemma \ref{lemmaupperface}, 
this orientation differs from the one of $J_{q^0 } \times J_{l_1 } \times \dots \times J_{l_{q_0}}$ exactly
by the sign $(-1)^{\sum_{j=1}^{q^0} (q^0-j)(l_j - 1) }$. As a consequence, 
the morphisms $\widetilde{\cal M} (H^1 \circ H^0)$ and ${\cal M} (H^0) \circ {\cal M} (H^1)$ coincide.
The result now follows along the same line as Theorem \ref{theocontinuation}. $\square$

\addcontentsline{toc}{part}{\hspace*{\indentation}Bibliography}

\bibliographystyle{abbrv}

\vspace{0.7cm}
\noindent 
Unit\'e de math\'ematiques pures et appliqu\'ees de l'\'Ecole normale sup\'erieure de Lyon,\\
 CNRS - Universit\'e de Lyon.

\end{document}